\documentclass[a4paper,12pt]{amsart}

\usepackage{gv}

\begin{document}

\title{Attractors of Piecewise Translation Maps}


\author{{D.~Volk}}
\email{dire.ulf@gmail.com}





\debug{\date{{\bf\today.}  File: {\jobname}.tex.}}

\begin{abstract}
Piecewise Translations is a class of dynamical systems which arises from some applications in computer science, machine learning, and electrical engineering. In dimension 1 it can also be viewed as a non-invertible generalization of Interval Exchange Transformations. These dynamical systems still possess some features of Interval Exchanges but the total volume is no longer preserved and allowed to decay.

Every Piecewise Translation has a well-defined attracting subset which is the locus of our interest. We prove some results about how fast the dynamics lock onto the attractor, geometry of the attractor, and its ergodic properties. Then we consider stochastic Piecewise Translations and prove that almost surely its attractor has zero Lebesgue measure. Finally we present some conjectures and supporting numerics.
\end{abstract}

\maketitle

\section{Introduction}  \label{s:intro}

Take a set $\Omega \subset \bbR^d$, cut $\Omega$ into finitely many pieces and independently move the pieces within $\Omega$ using only parallel translations of $\bbR^d$. The moved pieces are allowed to overlap thus the map needs not to be invertible. This is a \emph{piecewise translation (PWT)} map~$F$ defined on~$\Omega$.

Piecewise translations and a wider class of maps, piecewise isometries, have many applications in computer science, machine learning and electrical engineering: herding dynamics in Markov networks~\cite{Welling2009}, second order digital filters~\cite{chua1988chaos, deane2006piecewise}, sigma-delta modulators~\cite{ashwin2001dynamics, deane2002global, feely2000nonlinear}, buck converters~\cite{deane2006piecewise, deane2004buck}, three-capacitance models~\cite{Suzuki2004}, error diffusion algorithm in digital printing~\cite{Adler2005, Adler2012}.

In dimension 1 invertible PWTs are \emph{interval exchange transformations (IETs)} which are classic objects in ergodic theory but still attract a lot of interest. They have deep connections with polygon billiard maps, measurable foliations, translation flows, Abelian differentials, Teichm\"uller flows and other areas, see a review by Viana~\cite{viana2006ergodic}. Opposite to general PWTs, IETs preserve Lebesgue measure but their ergodic theory is still far from trivial.

Boshernitzan and Kornfeld in \cite{Boshernitzan1995} were first to consider general PWTs in dimension 1. They called them \emph{interval translation maps (ITMs)}. Like IETs, they are related to billiard maps if one allows one-sided straight mirrors scattered on the billiard table. Boshernitzan and Kornfeld showed that ITMs of rank less than 2 (i.e. that the endpoints and translation vectors of an ITM span a 2-dimensional subspace over the rational numbers) are finite type, that is their attractors stabilize in finite time. Such ITMs reduce to IETs after finitely many iterations. They also provided an example of a rank 3 ITM which was in fact infinite type.  They questioned the extent of which infinite type ITMs are typical.  This question was answered by Suzuki, Ito and Aihara in \cite{Suzuki2005} who showed that almost every double rotation (a piecewise rotation of two intervals, which are a class of ITMs of three intervals) are finite type, and later by Volk in \cite{Volk2014a} who demonstrated that almost every ITM of three intervals is finite type. This means that after finitely many iterates the dynamics converge on the attractor where the dynamics is a finite union of IETs. This allows to apply classic results such as unique ergodicity (Masur, Veech~\cite{masur1982interval, veech1982gauss}) and weak mixing (Avila, Forni~\cite{avila2007weak}) of almost every IET.

The study of piecewise isometries in 2 or more dimensions is still in its relative infancy, and most work was about piecewise isometries on the plane.
%
If a PWT in dimension $d>1$ has less than $d+1$ pieces, then $\bbR^d$ foliates into invariant affine spaces of dimension $d' < d$ and the PWT reduces to a disjoint collection of PWTs in dimension~$d'$. Thus the minimal interesting number of pieces is $d+1$.
In dimension 1, PWTs of $d+1=2$ intervals always stabilize after finitely many iterates. The attractor is a single closed interval. The dynamics on the attractor is an exchange of two subintervals, that is, a circle rotation.  In Section~\ref{s:pwt-in-r_d} we show that an analogous result holds in higher dimensions:
\emph{any} piecewise translation with $m=d+1$ pieces in $\R^d$ with rationally independent translation vectors is finite type, see Theorem~\ref{t:finite-type-erg}. The latter condition holds for almost every tuple of translation vectors. This theorem can also be seen as a direct generalization of a 2-dimensional result of Adler et al.~\cite{Adler2012}. See also the conjecture by Adler et al.~\cite{Adler2010} about the existence of an invariant fundamental set for the $n+1$ lattice.

\subsection{Plan of the paper}

In Section~\ref{s:pwt-in-r_d} we give all the necessary definitions and prove the first main result of this paper, Theorem~\ref{t:finite-type-erg}.

In Subsections~\ref{ss:geometry} and \ref{ss:partition} we generalize some of the results of Adler et al. to our settings. In Subsection~\ref{ss:geometry} we show that the Lebesgue measure of the attractor is an integer multiple of the volume of the torus factor. Note that in the setting of Adler et al~\cite{Adler2012}, this integer is shown to be equal to 1. This is what we observed in the numerical experiments in our case, too. The problem to prove this and some related conjectures (see Subsection~\ref{ss:geometry}) remains open and possibly related to the Pisot conjecture~\cite{akiyama2015pisot}.
\debug{\todo{say more..}}

Then in Subsection~\ref{ss:partition} we study ergodic properties of~$F$ of finite type. We show that for almost every~$x \in \Omega$ the frequencies of visits to partition elements are well defined and depend only on translation vectors. Moreover, they are equal to the normalized Lebesgue measures of the partition elements of the attractor.

Finally, in Subsection~\ref{ss:mgd1} we discuss the case when a PWT in~$\bbR^d$ has the number of partition elements $m$ bigger than $d+1$. In dimension 1, it is known since the pioneering work of Boshernitzan and Kornfeld~\cite{Boshernitzan1995} that for $m=3$ the finite and infinite type regimes coexist in the same parameter space. Volk~\cite{Volk2014a} showed that for $m=3$ all the finite type ITMs are contained in an open and dense subset of full measure of the parameter space, see also~\cite{Bruin2012}. This complements the explicit construction of a continuum of parameter values corresponding to ITMs of infinite type~\cite{Boshernitzan1995}.
Similar results also hold for certain families with $d=1, m > 3$, see~\cite{Schmeling2000}.
\debug{\todo{Add 1--2 more refs}}

Our numerical experiments suggest that the situation in $d=2, m = 4$ is similar to the case $d=1, m=3$. Namely, we observed that for fixed partitions, the parameter space of admissible translation vectors consists of big regions of finite type with relatively small stabilization time, and narrow regions where the stabilization time explodes. This suggests the parameters are near infinite types. In Subsection~\ref{ss:mgd1} we present some pictures in support of this conjecture.

\debug{
Later, in a preprint~\cite{Adler2012} \\
\textbf{Result A (Ergodic Inputs).} For acute simplices the minimal absorbing invariant set for the error diffusion with an ergodic constant input is a fundamental set for the lattice generated by the simplex.
\\
We propose but cannot yet prove the following
\textbf{Conjecture.} The ergodicity assumption in Theorem I.16 is redundant.
\\
\textbf{Result B (Sub-Tiles).} If a bounded forward invariant set of a generalized (arbitrary partition) error diffusion on a simplex is fundamental for the simplex lattice, then each part of this invariant set (i.e., each intersection with the partition) is a
fundamental domain for a derived lattice. (pre-requisites are ``we only know these assumptions to hold true in the acute case when $d = 1, 2$ and in the acute case with Ergodic Input when $d > 2$.''

Is this true in our setting?
}

\section{Piecewise Translation Maps in $\mathbb{R}^d$}  \label{s:pwt-in-r_d}

\subsection{Definitions and main players} \label{ss:defs}

We say a compact set in affine space~$\mathbb{R}^d$ is a \emph{region} if it is the closure of its interior. Consider a region
$\Om \subset \mathbb{R}^d$ and its partition~$\mathfrak{P}$ into $m$ smaller regions $\Om = P_0 \cup \dots \cup P_{m-1}$ such that the pieces intersect only at the boundaries. Also assume that $\Leb \partial P_i = 0$ for all $P_i$.

Now to each~$P_i$ we attach a vector~$v_i\in \R^d$ such that for any $x\in P_i$ we have $x+v_i \in \Om$. For $x \in P_i$ we denote $i(x) = i$, $P(x) = P_{i(x)}$ and $v(x) = v_{i(x)}$. Note there is some ambiguity at overlapping boundaries of partition elements. But because overlaps have zero measure, we just assume this ambiguity is resolved in some fixed way, so that $i(x)$, $P(x)$ and $v(x)$ are well-defined measurable functions.
Then the map $F\colon\Om\to \Om$, $F(x)=x+v(x)$, is a \emph{piecewise translation of $m$ branches}, and $\mathfrak{v} = (v_0, \dots, v_{m-1})$ are its \emph{translation vectors}.

\begin{figure}[h]
\begin{center}
\includegraphics[width=2.5in]{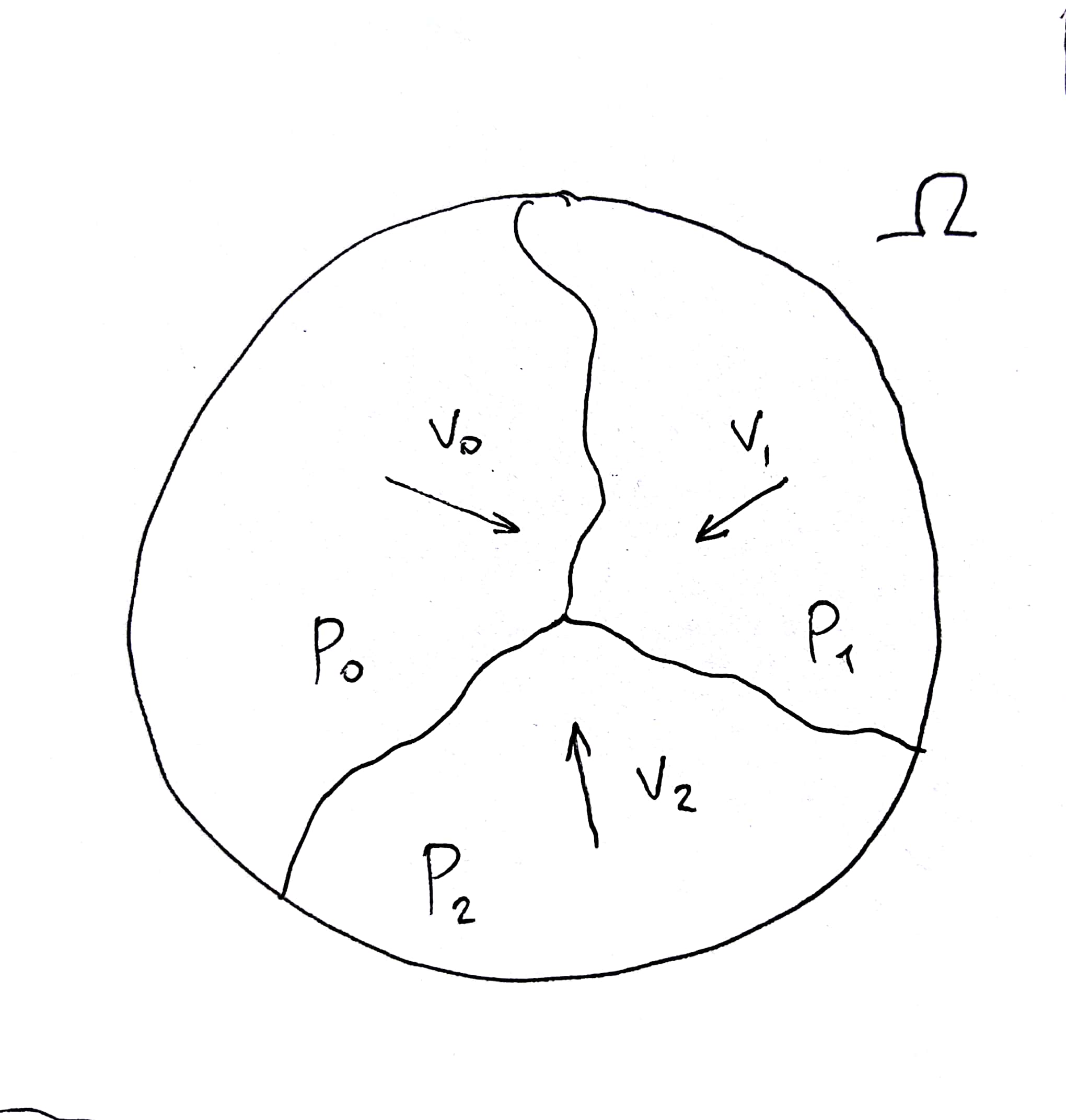}
\end{center}
\caption{Piecewise Translation Map of a disk in~$\bbR^2$ with 3 branches.}
\label{f:pwt}
\end{figure}

Denote $\Om_0 = \Om$ and let $\Om_n = \overline{F(\Om_{n-1})}$ for $n\geq 1$. Then~$\Om_0 \supseteq \Om_1 \supseteq \Om_2 \dots$ is a sequence of nested non-empty compact sets.
By Cantor's intersection theorem, this sequence has non-empty compact intersection~$A$, which we call the \emph{attractor} of $F$. The set~$A$ is $F$-invariant. Note that~$A$ is independent of a particular way of resolution of the ambiguity about~$i(x)$. We say a piecewise translation is \emph{finite type} if the sequence $(\Om_n)$ stabilizes, i.e., $\Om_n=\Om_{n+1}$ for some~$n\in\N$.  Otherwise we shall say it is \emph{infinite type}. This terminology mirrors the one proposed by Boshernitzan and Kornfeld~\cite{Boshernitzan1995} in the context of interval translation maps. Obviously, attractors of PWTs of finite type are regions. Schmeling and Troubetzkoy~\cite{Schmeling2000} proved that in dimension 1, transitive attractors of infinite type are Cantor sets.

Our goal in this Section is to study attractors of piecewise translations.
We will be most interested in the simplest non-trivial case, $m = d+1$ and $\rank\mathfrak{v} = d$. This case is also special because it admits a \emph{torus rotation factor} first introduced by Adler et al. in~\cite{Adler2010} for a special class of piecewise translations.

Namely, let $\mathcal{L}$ be the lattice isomorphic to~$\Z^d \subset \R^d$ generated by the vectors $(v_1-v_0),(v_2-v_0),\dots,(v_d-v_0)$, and let $\T$ be the $d$-dimensional torus~$\R^d/\mathcal{L}$. Take $\pi \colon \R^d \to \T$ to be the canonical projection.  The crucial observation is that $\pi$ is a semiconjugacy of $F$ and the torus rotation map $R\colon\T\to\T$,  $R(\phi) =\phi+ v_0$, i.e., $R\circ \pi=\pi\circ F$.
Indeed,
$$
\pi\circ F(x)= \pi(x+v(x)) = \pi(x +(v_i-v_0)+v_0) = \pi(x+v_0)=\pi(x) +v_0.
$$
We will usually assume that vectors~$\mathfrak{v}$
are rationally independent, i.e., one has $n_1 v_1 + \dots + n_m v_m = 0$ for $n_i \in \bbZ$ iff all $n_i = 0$.
This is equivalent for~$R$ to be ergodic. This assumption holds on a full measure subset of the space of vector parameters.

\subsection{Finiteness theorem for $m = d+1$}

In the realm of piecewise translations people are often interested whether they are dealing with maps of finite or infinite type. For certain parameter values this question is easy to answer. For instance, if $\mathfrak{v}$ are rational vectors, then every orbit is finite, of uniformly bounded length. Thus $F$ is finite type.

On the other hand, it is a straightforward exercise to show that in dimension~1, \emph{any} piecewise translation of 2~branches is finite type. In this case, the attractor~$A$ is an interval, and $F|_A$ is an exchange of two intervals which is a circle rotation. Theorem~\ref{t:finite-type-erg} is a far-going generalization of this fact.

From now and until Subsection~\ref{ss:mgd1} we always assume that $F$ is a piecewise translation in~$\bbR^d$ with $m = d+1$ branches.

\begin{Thm}	\label{t:finite-type-erg}
If~$R$ is ergodic, then $F$ is finite type.
\end{Thm}

\begin{figure}[h]
\begin{center}$
\begin{array}{cc}
\includegraphics[width=2.5in]{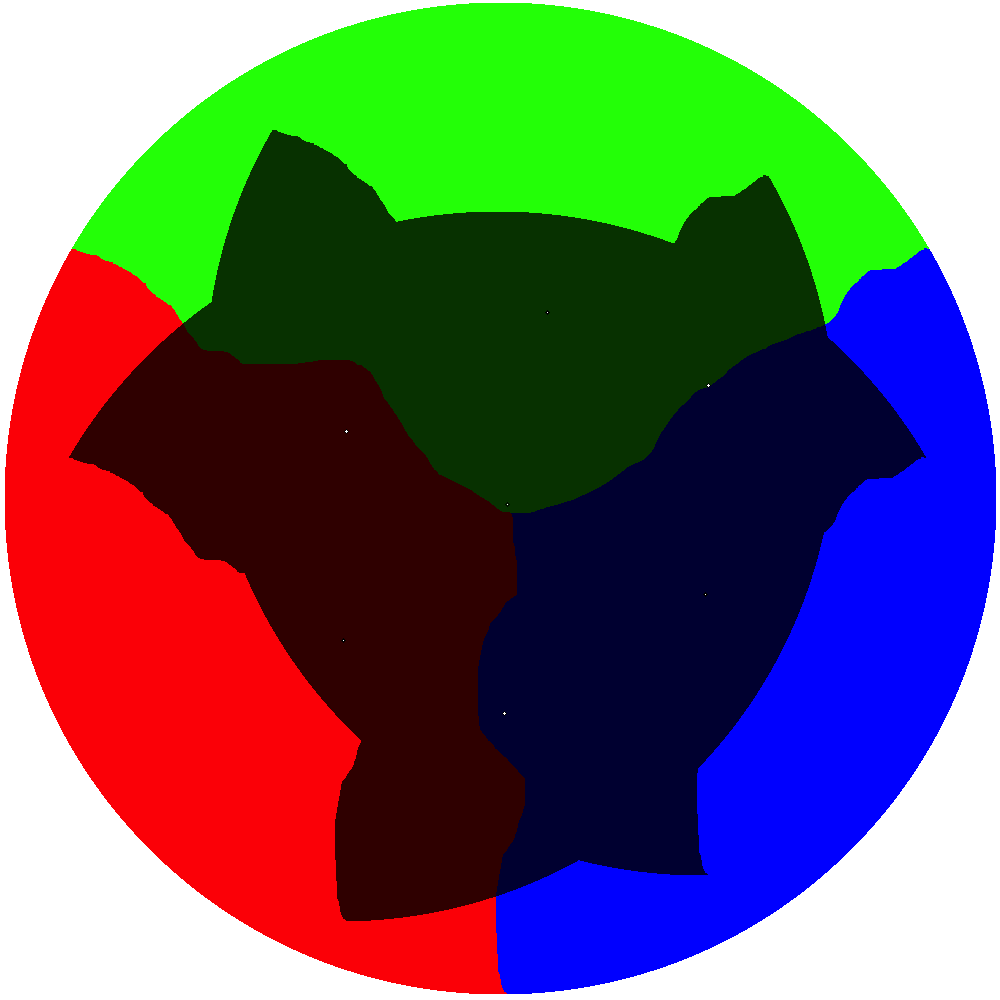} &
\includegraphics[width=2.5in]{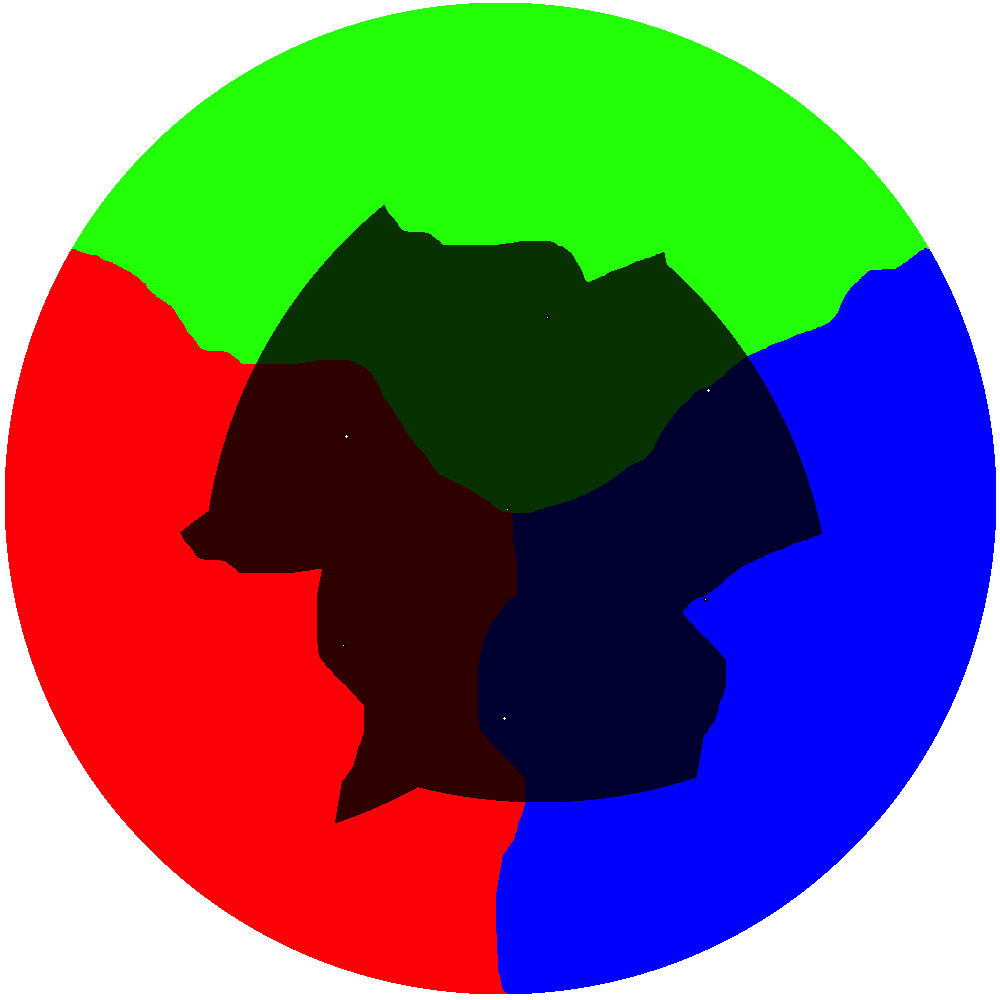} \\
\includegraphics[width=2.5in]{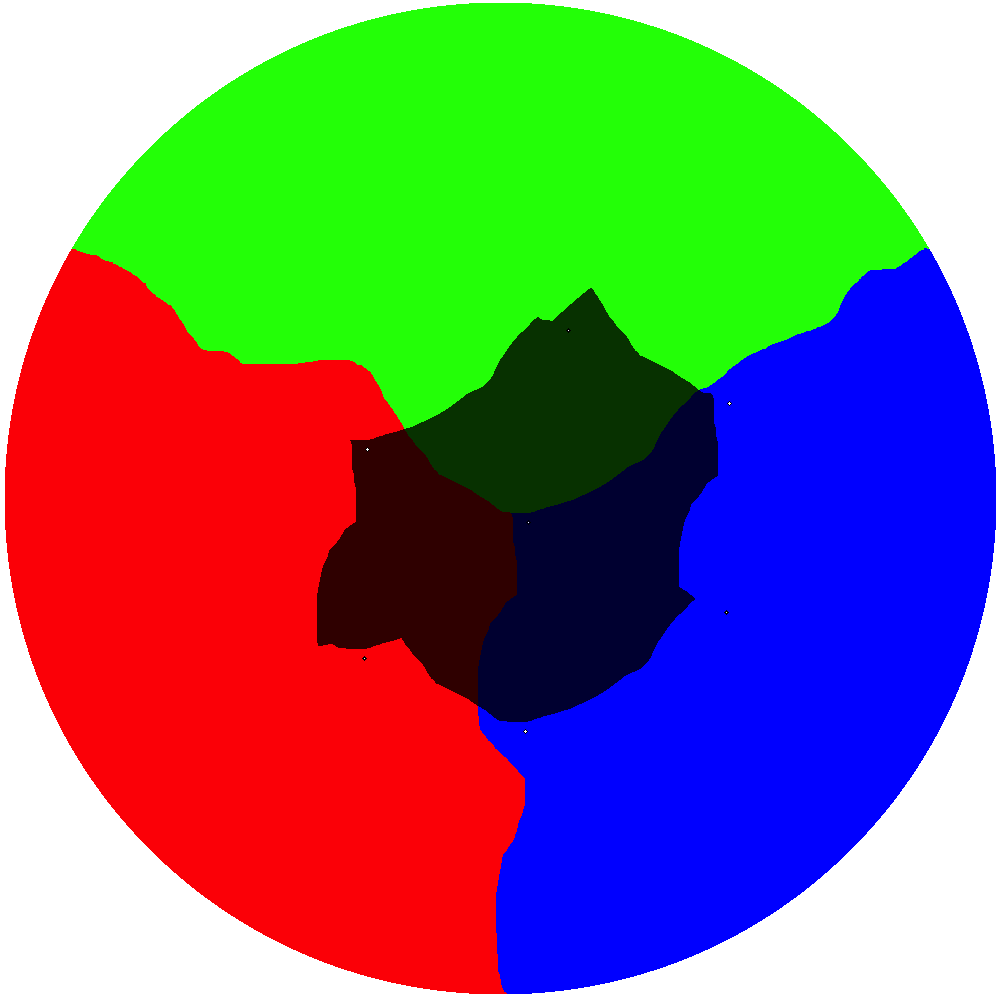} &
\includegraphics[width=2.5in]{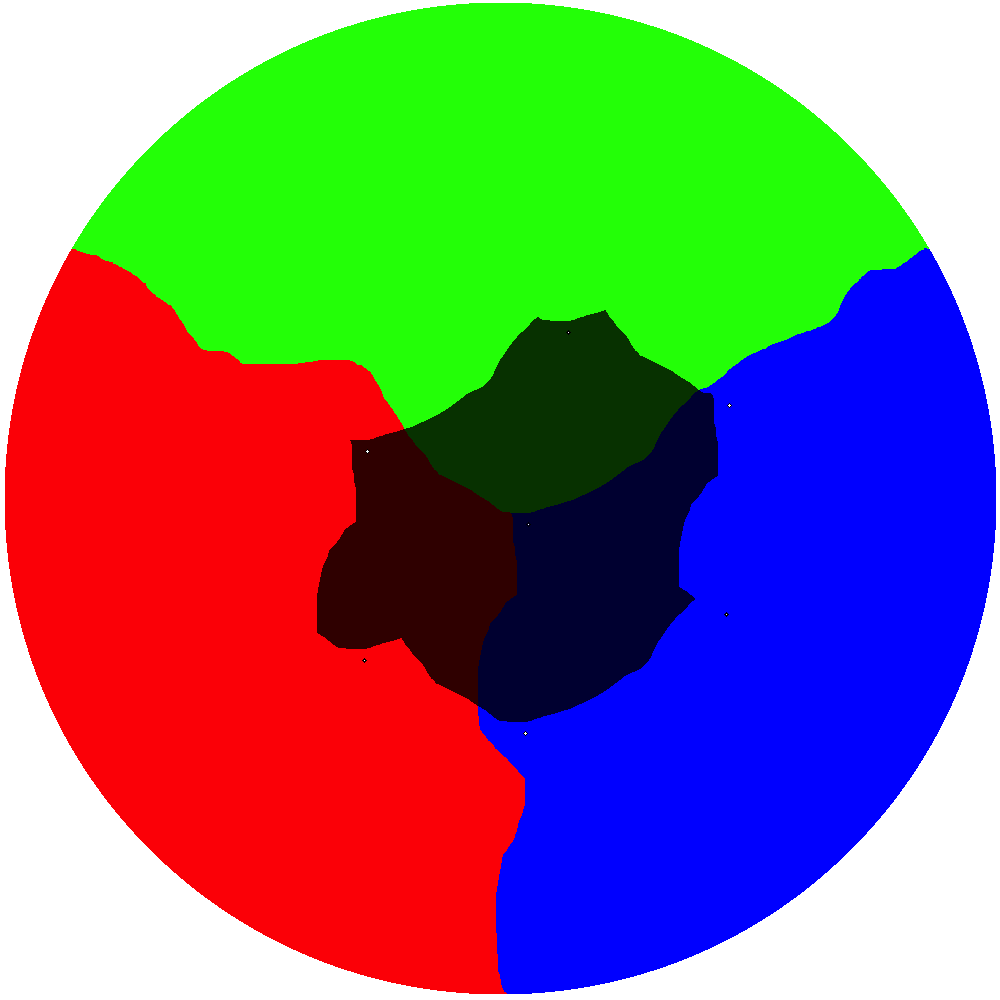}
\end{array}$
\end{center}
\caption{Finite type attractor for $m = 3$ in $\bbR^2$.\\ The disk~$\Om$ converges to attractor in 5 iterates. $F^n(\Om)$ for $n = 1,2,4,5$.}
\label{f:disk-3}
\end{figure}

\begin{proof}
First, let us show that $\pi(A)=\T$. Because $A$ is compact and non-empty and $\pi$ is continuous, $\pi(A)$ is compact and non-empty.  Since $A$ is invariant under $F$, $\pi(A)$ must be invariant under $R$. The ergodic rotation $R$ has only two invariant closed sets, the whole torus and the empty set. Therefore $\pi(A)=\T$.

On the other hand, $\partial A$ and $\pi (\partial A)$ are not invariant.
By definition, $\partial A$ is closed nowhere dense in~$\bbR^d$. Since $\Omega$ is compact, the projection $\pi\colon\Omega\to\bbT$ is at most finitely branched, so locally $\pi(\partial A)$ is a finite union of nowhere dense sets, and thus is nowhere dense in~$\bbT$, too. Hence we can take a small ball~$B$ in the complement to $\pi (\partial A)$, and $\varepsilon > 0$ such that for the $\varepsilon$-neighborhood $U \subset \T$ of $\pi(\partial A)$ we have $U \cup B = \emptyset$.

Let $V$ be the $\varepsilon$-neighborhood of $A$ upstairs, and $V' \subset V$ be the $\varepsilon$-neighborhood of $\partial A$. By definition of~$A$ there exists $N > 0$ such that $\forall n > N$ we have $\Om_n \subset V$. Thus $\forall n > N$ $\forall x \in \Omega$ $F^n(x) \in V$.

It is well-known that any ergodic torus rotation is uniformly minimal, that is, for any open~$B \subset \T$ there exists $M>0$ such that $\forall \phi \in \T$ we have $R^k \phi \in B$ for some $0 < k < M$.

In particular, we can take $\phi = \pi (F^n(x))$. Then for $y = F^{n+k}(x)$ we know that $y \in V$, and $y \notin V'$. Thus $y \in \inter A$. So $\exists N,M$ such that $\forall x \in \Om$ $F^{n+k} (x) \in \inter A \subset A$ for some $n+k < N+M$. Therefore $F$ is finite type.
\end{proof}

\begin{Cor}
If~$R$ is ergodic, then $A$ is a region, and $F|_A$ is an exchange of $m+1$ pieces semiconjugated to $R$.
\end{Cor}


To us, the Theorem and its Corollary were wonderful and hard to believe at first, because in the dimensions $d>1$, it is rather difficult to construct nontrivial region exchange maps at all, and here you get a region exchange map out of nowhere. Moreover, by construction, the boundary~$\partial A$ is a finite union of translated copies of some parts of~$\partial P_i$.

\begin{Pbm}
Is it possible to identify these pieces \emph{a priori}, before knowing the attractor itself?
\end{Pbm}

\subsection{Geometry of attractors in $d = m+1$}	\label{ss:geometry}


For any~$\phi \in \T$, define $\xi (\phi) = \#\{\pi^{-1}(\phi)\cap A\} \in \bbN$.
\begin{Thm}\label{t:finite-cover}
If~$R$ is ergodic, then $\xi(\phi) = const$ for almost all $\phi$.
\end{Thm}
\begin{proof}
Denote by $J \subset \T$ the subset of points whose orbits never come to $\pi(\partial A)$. Note that $\Leb J = \Leb \T$. Because the map~$F$ is 1-1 on $\pi^{-1} J \cap A$, the function~$\xi$ is invariant under $R$ almost everywhere, thus constant.
\end{proof}

\begin{Cor}\label{c:Avol}
If~$R$ is ergodic, then the volume of $A$ is an integral multiple of the volume of $\T$, i.e. $\frac{\Leb A}{\Leb \T} = \l \in \bbN$.
\end{Cor}

Now the natural question is, what are the possible values of $l$. Our numerical experiments suggest that $l$ always equals to 1. Another numerical observation is that $A$ is always a tile: one can cover the whole~$\bbR^d$ with the translated copies~$A+\mathcal{L}$ with overlaps only at boundaries. The same holds true for every $A_i = A \cap P_i$ (for different lattices).

\begin{figure}[h]
\begin{center}
\includegraphics[width=5in]{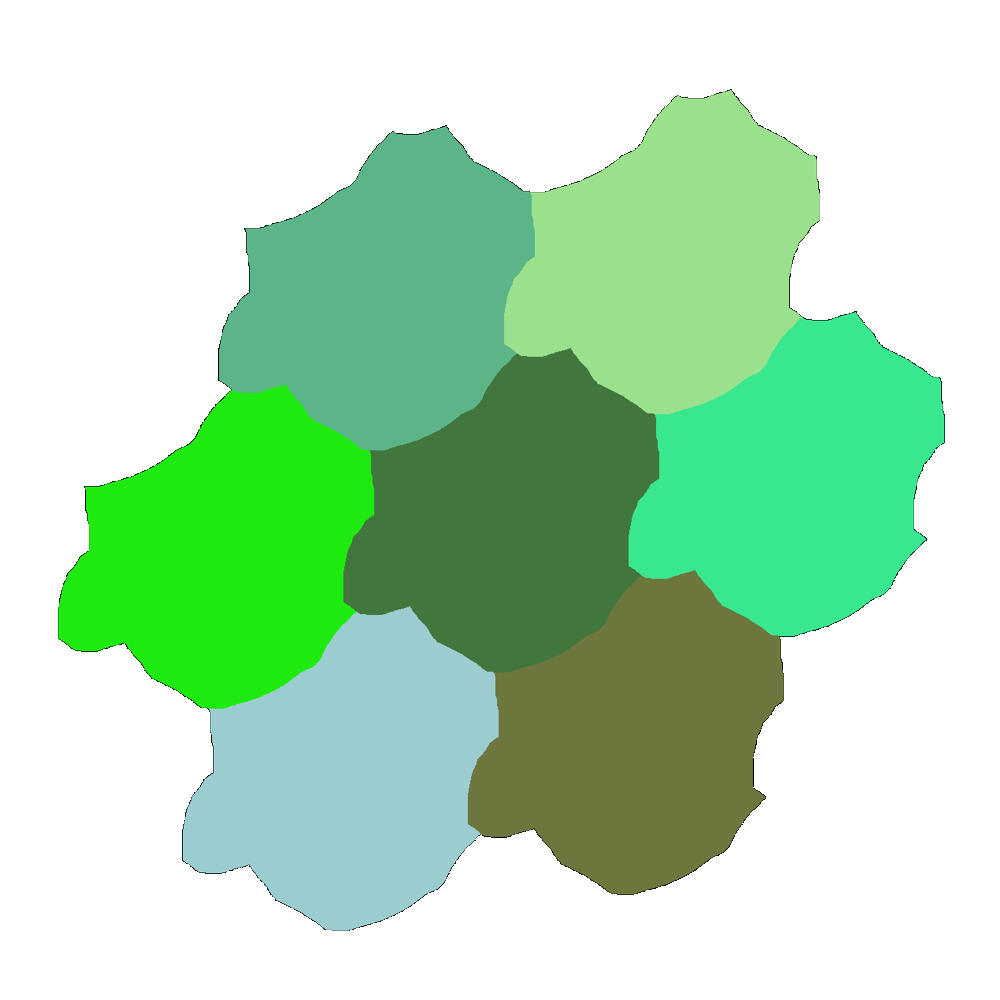}
\end{center}
\caption{The attractor as a tile}
\label{f:tile}
\end{figure}

Adler et al.~\cite{Adler2012} proved these statements for a very special class of partitions and translation vectors in~$\bbR^2$. Namely, their partitions are Voronoi domains for three points $p_0,p_1,p_2$ in the plane which form a non-obtuse triangle, and their translation vectors are $\gamma-p_i$ for a $\gamma$ inside the triangle $(p_0, p_1, p_2)$. However, the general question remains open.

\subsection{The partition of the attractor for $m = d+1$}	\label{ss:partition}

In this Subsection we do not assume that~$R$ is ergodic.

\begin{Prop}	\label{p:inf-times}
For every~$x \in \Om$, its orbit~$F^n(x)$ visits every partition element~$P_i$ infinitely many times.
\end{Prop}

\begin{proof}
Assume that for some $P_j$ there exists~$N$ such that for any $n > N$ we have $F^n(x) \notin P_j$. On the other hand, at least one~$P_k$ \emph{must} be visited infinitely many times.

Because~$\rank \mathfrak{v} = d$, every $d$ vectors from~$\mathfrak{v}$ are a basis in~$\bbR^d$. Take~$\mathfrak{v} \setminus \{ v_j \}$ as the basis.

In this basis, consider the coordinates of the orbit~$F^n(x)$. Note that when~$n > N$,
at every iterate some coordinate increases by 1, and the others stay the same. Moreover, the $k$-th coordinate tends to~$+\infty$.
On the other hand, the assumption~$F(\Omega) \subset \Omega$ implies that every orbit is bounded. This contradiction proves the proposition.
\end{proof}

Because~$\mathfrak{v} = \{v_0, \dots, v_d\}$ is a system of $d+1$ vector of rank~$d$, there exists a unique (up to a factor) nontrivial tuple~$(\alpha_0, \dots, \alpha_d)$ such that $\sum \alpha_i v_i = 0$. Again because $\rank \mathfrak{v}  = d$, we have~$\alpha_i \ne 0$ for all~$i$.

\begin{Prop}	\label{p:alpha-i}
One can choose the factor so that $\alpha_i > 0$ for every~$i$, and $\sum \alpha_i = 1$.
\end{Prop}

\begin{proof}
Assume the contrary. Then there exists a tuple~$(\alpha_0, \dots, \alpha_d)$ such that $\sum \alpha_i v_i = 0$, and $\alpha_j < 0$, $\alpha_k > 0$ for some $j \ne k$. Because~$\rank \mathfrak{v} = d$, every $d$ vectors from~$\mathfrak{v}$ are a basis in~$\bbR^d$. Take~$\mathfrak{v} \setminus \{ v_j \}$ as the basis.

In this basis, let us look at the $k$-th coordinate of all the vectors. First, the $k$-th coordinate of all but~$v_j$, $v_k$ equals zero. Because $v_j = \frac{1}{-\alpha_j} \sum_{i \ne j} \alpha_i v_i$, the $k$-th coordinate of $v_j$ is $\frac{\alpha_k}{-\alpha_j} > 0$. Finally, the $k$-th coordinate of $v_k$ is 1. In particular, every~$v_i$ has its $k$-th coordinate bigger or equal to zero.

By Proposition~\ref{p:inf-times}, for every~$x \in \Om$ the orbit~$F^n(x)$ visits every partition element~$P_i$ infinitely many times. Then its $k$-th coordinate grows arbitrary large and thus $F^n(x)$ escapes to infinity. But we know that every orbit is bounded. This contradiction proves the proposition.
\end{proof}

The following lemma tells us that for $F$ of finite type, the coefficients~$\alpha_i$ have a nice geometric meaning.
\begin{Lem}	\label{l:alpha}
Let $m=d+1$, and assume~$F$ is finite type. Then
\begin{enumerate}
\item\label{i:alpha-freq} for Lebesgue-almost every $x\in \Omega$ the frequency of visits of $x$ to the partition element~$P_i$ is well defined and equals $\alpha_i$;
\item\label{i:alpha-area} the normalized Lebesgue measures of attractor pieces~$\frac{\Leb A_i}{\Leb A}$ equal $\alpha_i$.
\end{enumerate}
\end{Lem}

This result can be viewed as a step towards proving that $F|_A$ is ergodic provided $R$ is ergodic.

\begin{proof}
Because~$F$ is finite type, we can assume from the beginning that~$x \in A$ and that~$\Leb A = 1$.
Let~$b_i(x)$ be the indicator functions of $A_i$, and $b_i^{(k)}(x) = \sum_{n=0}^{k-1} b_i(F^n(x))$ be their sums along the orbit of $x$. Note that the $k$-th iterate of $x$ has the following form:
$$
F^k (x) = x + b^{(k)}_0(x) v_0 + \dots + b^{(k)}_d(x) v_d, \text{ and } \sum_{i=0}^d b^{(k)}_i(x)  = k.
$$

As the Lebesgue measure is invariant under~$F|_A$, we can take its ergodic decomposition~$\Leb = \int \mu_\eta \,d\eta$. Then for Lebesgue-almost every point~$x \in A$ there exists an ergodic measure~$\mu_\eta$ such that $x$ is $\mu_\eta$-generic, and thus the Birkhoff averages $\frac1k b_i^{(k)}(x)$ converge. Let~$\beta_i(x) = \lim_{k\to+\infty} \frac1k b_i^{(k)}(x)$. By construction, $\beta_i(x)$ is the frequency of $x$'s visits to $P_i$, and $\sum_i \beta_i(x) = 1$.

Now consider the limit
$$
\lim_{k\to+\infty} \frac1k F^k(x) = \lim_{k\to+\infty} \frac1k (x + b^{(k)}_0(x) v_0 + \dots + b^{(k)}_d(x) v_d) =
\beta_0(x) v_0 + \dots + \beta_d(x) v_d.
$$
Because the orbit~$F^k (x)$ stays in~$\Omega$ and thus is bounded, $\frac1k F^k (x) \to 0$ as $k \to \infty$. Thus $\sum \beta_i(x) v_i = 0$, and by construction~$\sum \beta_i(x) = 1$. Such a tuple is unique, therefore $\beta_i(x) = \alpha_i$ for all~$i$ and all~$x \in A$ (and all~$x \in \Omega$ too).

On the other hand, for every~$\mu_\eta$ take a $\mu_\eta$-generic $x$. By Birkhoff ergodic theorem, $\mu_\eta(A_i) = \beta_i(x) = \alpha_i$, and thus $\Leb(A_i) = \int \mu_\eta(A_i) \,d\eta = \alpha_i$.
\end{proof}

\subsection{Fate maps and Diagrams}

\subsubsection{Fate map}

Given a piecewise translation map~$F$, we have a partition~$\mathfrak{P}$ from the beginning. The fate map $\mathcal{F} \colon \Om \to \Sigma$, where $\Sigma = \{0,\dots,d \}^{\bbN_0}$ is the one-sided symbolic space associated with the partition~$\mathfrak{P} = (P_0, \dots, P_d)$, and $\bbN_0 = \bbN \cup \{0\}$, is given by
$$
\mathcal{F}  \colon x \mapsto (i(x),i(F(x)),i(F^2(x)),\dots).
$$
By construction, for the left shift~$\sigma \colon \Sigma  \to \Sigma $ we have $\mathcal{F}  \circ f = \om \circ \mathcal{F}$.


Lemmas~\ref{l:fate-aper} and~\ref{l:fate-per} hold true for any number~$m$ of partition elements.

\begin{Lem}	\label{l:fate-aper}
For any aperiodic $\om \in \Sigma$ we have $\Leb(\mathcal{F}^{-1} (\om)) = 0$.
\end{Lem}

\begin{proof}
Note that $\mcF^{-1} (\om) \stackrel{f}{\longrightarrow} \mcF^{-1} (\sigma\om)$ is a 1-1 rigid translation, because whole $\mcF^{-1} (\om)$ lies in the same element of the partition. Thus for any $n_1 < n_2 \in\bbN_0$ there exists $v \in \bbR^n$ such that 
$$
F^{n_1}(\mcF^{-1} (\om)) =t v + F^{n_2}(\mcF^{-1} (\om)).
$$
In particular, 
$$\
\Leb(F^{n_1}(\mcF^{-1} (\om))) = \Leb(F^{n_2}(\mcF^{-1} (\om))).
$$
By the definition of the fate map, $\forall \om_1 \ne \om_2$ we have $\mcF^{-1} (\om_1) \cap \mcF^{-1} (\om_2) = \emptyset$. For an aperiodic $\om$ there exists a sequence of numbers $n_1 < \dots < n_k \dots \in \bbN_0$ such that all $\sigma^{n_k}\om$ are different. Thus $F^{n_k} (\mcF^{-1} (\om))$ is countably many pairwise disjoint sets of same Lebesgue measure, confined in a region~$\Om$ of finite volume. This contradiction proves the lemma.
\end{proof}

\begin{Lem}\label{l:fate-per}
If any point $x \in \Om$ has a periodic fate, then the translation vectors~$\mathfrak{v}$ must be rationally dependent.
\end{Lem}

\begin{proof}
Assume for some $x \in \Om$ we have $\mcF(x) = \om$, and $\sigma^p\om = \om$.
The map $\mcF^{-1}(\om) \stackrel{F^p}{\longrightarrow} \mcF^{-1}(\om)$ is 1-1, moreover, it is a parallel translation by a nontrivial combination of vectors~$v_i$, $m_0 v_0 + \dots + m_d v_d$, $m_i \in \bbN_0$. Because~$F^p$ leaves the whole set~$\mcF^{-1}(\om)$ fixed, we have $m_0 v_0 + \dots + m_d v_d = 0$. Moreover, because $\rank \mathfrak{v} = d$, we have $m_k > 0$ $\forall k$.
\end{proof}

%

\begin{Lem}
Let $m=d+1$, and assume~$F|_A$ is minimal. Then the fate map~$\mathcal{F}$ is injective.
\end{Lem}

\begin{proof}
Let~$x \ne y$, $x,y \in \Om^*$, have the same fate~$\om$. This means that every next iterate is the \emph{same} translation for orbits of $x$ and $y$. Thus for any~$n \ge 0$ we have $F^ny - F^nx = y - x$.
By Theorem~\ref{t:finite-type-erg}, there exists~$M$ be such that for any~$n \ge M$ the points~$F^n x$ and $F^n y$ belong to the attractor~$A$.
Let~$z_0 = \argmax_{z \in A} (y-x, z)$ where $(.,.)$ is the standard scalar product. Let~$U$ be the $\frac12 |y-x|$-neighborhood of~$z_0$.
Because~$A$ is a region, we can take a small ball~$B \subset \inter A \cap U$.

\noindent By minimality of~$F|_A$, there exists $n > M$ such that $F^n x \in B$. But then $F^n y = F^n x + (y-x)$ is outside of~$A$! This controversy proves the lemma.
\end{proof}

%
%

\subsubsection{Diagrams}

Let $x \in A$. We define the \emph{diagram} of $x$ to be the subset~$D(x) \subset \mathcal{L}$ such that $\lambda \in D(x) \iff x + \lambda \in A$. In other words, $D(x) = \pi^{-1} (\pi(x)) - x$. By Theorem~\ref{t:finite-cover}, $\# D(x) = l$ for every $x \in A$. The \emph{class} of $x$ is the subset~$C(x) \subset A$ of the points with the same diagram.

A \emph{free diagram} is a class of equivalence of diagrams under parallel translations. The \emph{free class} of $x$ is the subset~$\hat C(x) \subset A$ of the points with the same free diagram~$\hat D$. Obviously, $\hat C(x) = \{ y \in A \,|\, y-x \in \mathcal{L} \}$. In the following proposition we summarize some basic properties of diagrams and classes.

\begin{Prop}
Let $F$ be a piecewise translation with $m = d+1$ branches, and assume~$R$ is ergodic. Then
\begin{enumerate}
\item $\forall x \in A$ $C(x)$ is a region.
\item if $D(x) \ne D(y)$, then $C(x) \cap C(y) = \emptyset$.
\item if $\hat D(x) = \hat D(y)$, then $C(y) = C(x) + \lambda$ for some $\lambda \in \mathcal{L}$.
\item $\hat C(x) = \{ y \in A \,|\, y-x \in \mathcal{L}\}$.
\item $\hat C(x) = \bigcup\limits_{\hat D(x) = const} C(x)$. Note that this a union of $l$ identical translated pieces which overlap at most at the boundaries.
\end{enumerate}
\end{Prop}

\subsection{More than $d+1$ pieces in~$\bbR^d$}	\label{ss:mgd1}

For piecewise translation maps with $m$ branches in~$\bbR^d$, where $m > d+1$, we have more questions than answers.

\subsubsection{Finiteness problem}

To begin with, we know from Boshernitzan and Kornfeld~\cite{Boshernitzan1995} that in the dimension 1 there exist piecewise translations of 3 branches with a Cantor attractor. Their Cartesian product gives a piecewise translation of 9 branches on the plane, with the product Cantor attractor, too. It is unknown if this is possible with a fewer number of branches.

\begin{Pbm}
Give at least one explicit example of a Cantor attractor for $m = 4$ in $\bbR^2$.
\end{Pbm}

\begin{figure}[h]
\begin{center}$
\begin{array}{cc}
\includegraphics[width=2.5in]{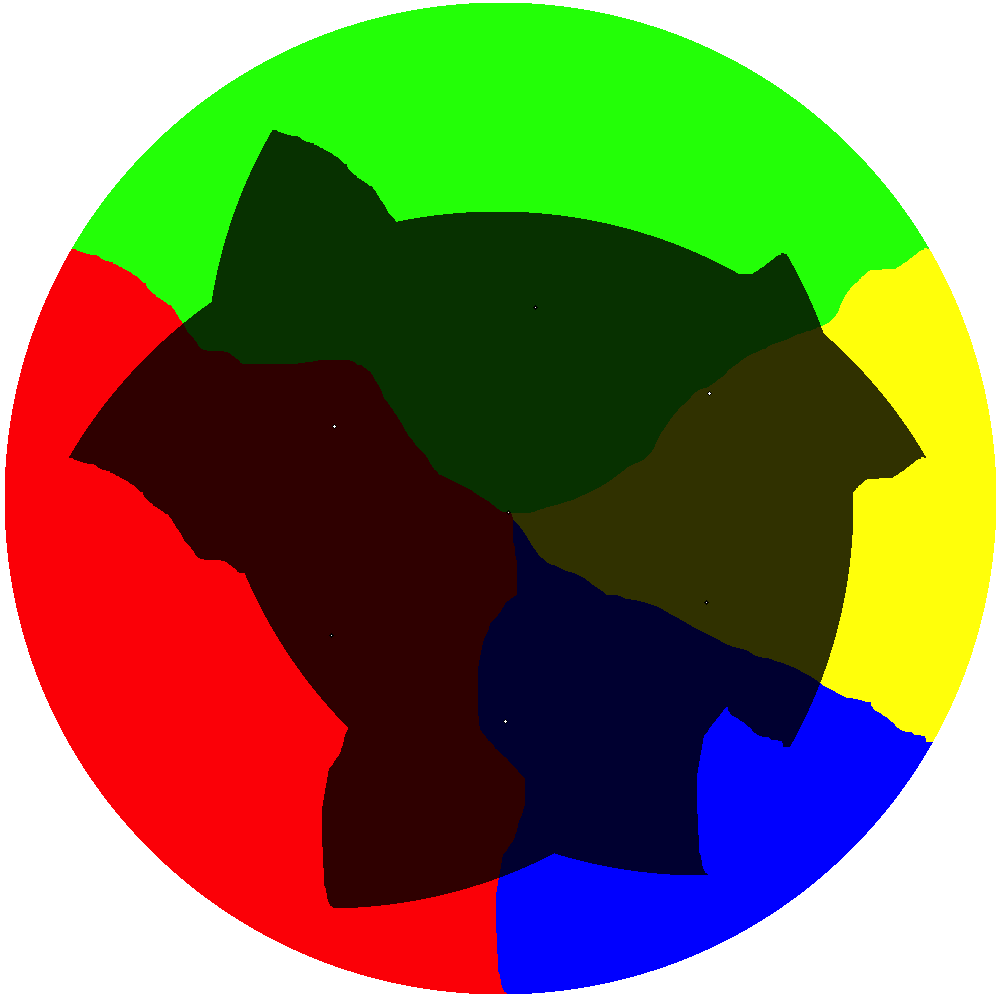} &
\includegraphics[width=2.5in]{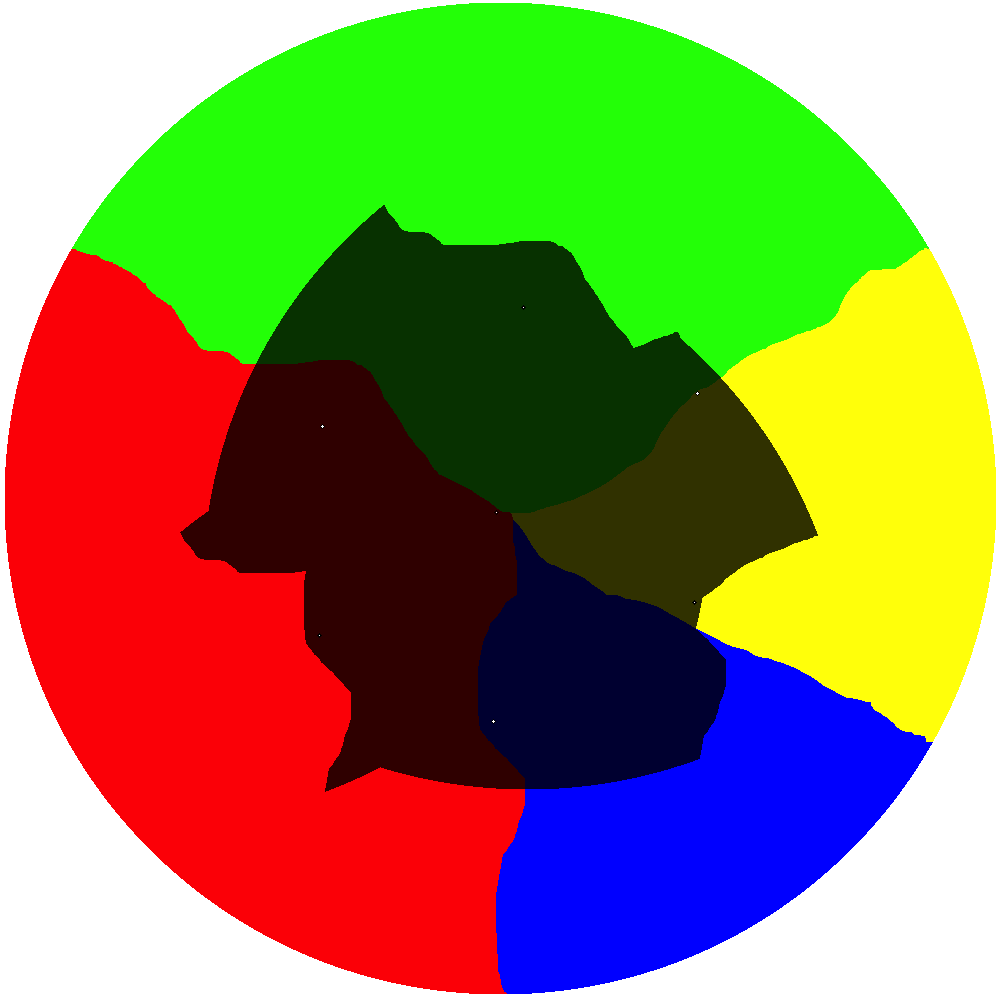} \\
\includegraphics[width=2.5in]{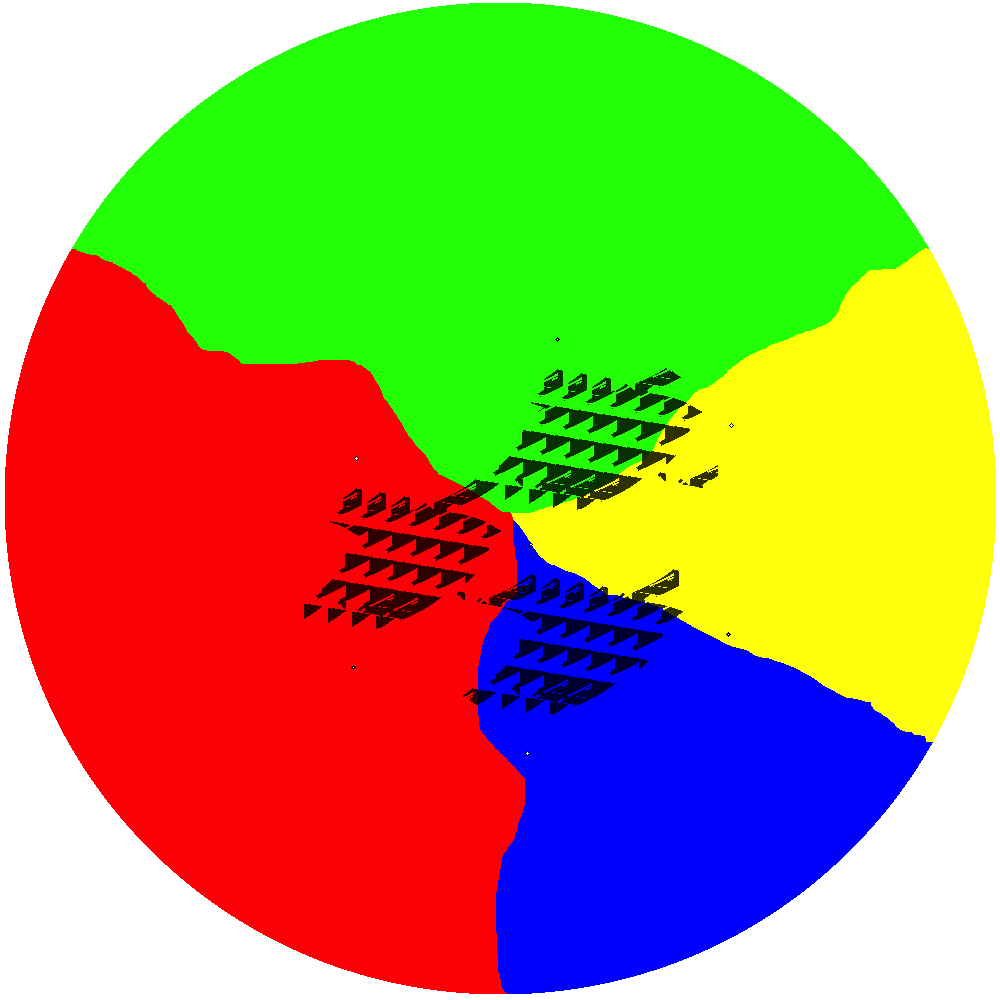} &
\includegraphics[width=2.5in]{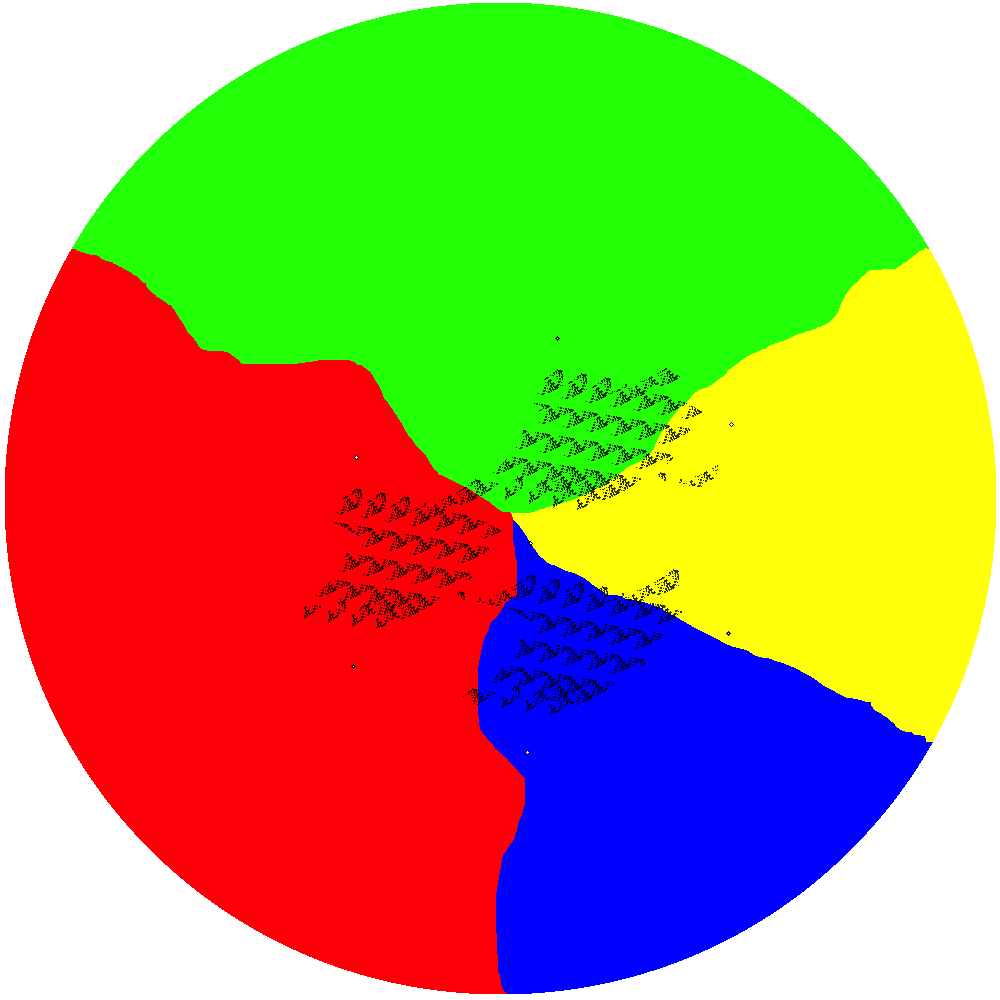}
\end{array}$
\end{center}
\caption{Seemingly Cantor-like attractor for $r = 4$ in $\bbR^2$.\\
The disk converges to attractor (numerically) in 2581 iterates.\\
$F^n (\Om)$ for $n = 1, 2, 100, 2581$.}
\label{f:disk-4-inf}
\end{figure}

In~\cite{Volk2014a}, Volk proved that in the dimension 1, the set of 3-branched piecewise translations of infinite type has zero measure in the parameter space.

\begin{Pbm}
Fix $m > d+1$ and a partition~$\Om = P_0 \cup \dots \cup P_{m-1}\subset \bbR^d$. Consider the set~$\mathcal{V} = \{\mathfrak{v}\}$ of all possible tuples of translation vectors which make it a well-defined piecewise translation. This is a region in~$(\bbR^d)^m$. Endow~$\mathcal{V}$ with a Lebesgue measure. Is it true that
\begin{itemize}
\item almost all piecewise translations are finite type?
\item the subset of infinite type piecewise translations has Hausdorff dimension less than $d \cdot m$?
\end{itemize}
\end{Pbm}

\begin{figure}[h]
\begin{center}$
\begin{array}{cc}
\includegraphics[width=2.5in]{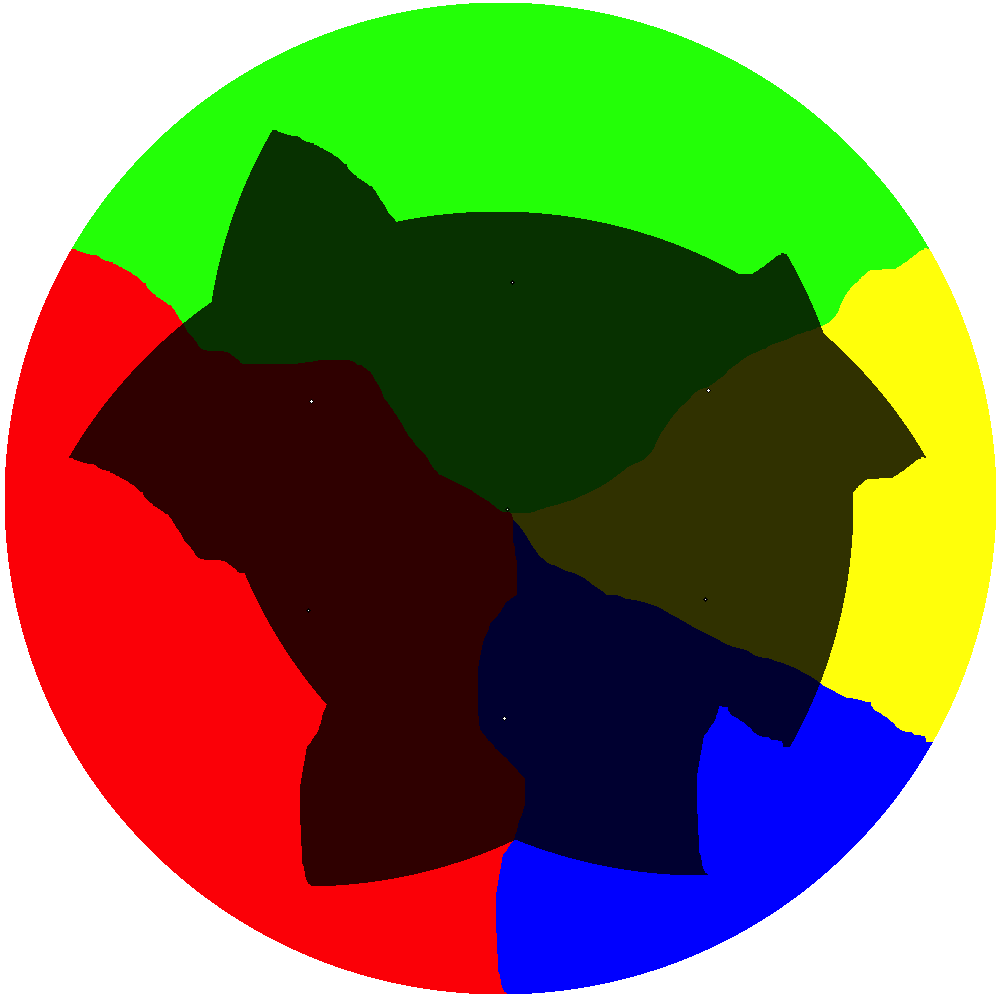} &
\includegraphics[width=2.5in]{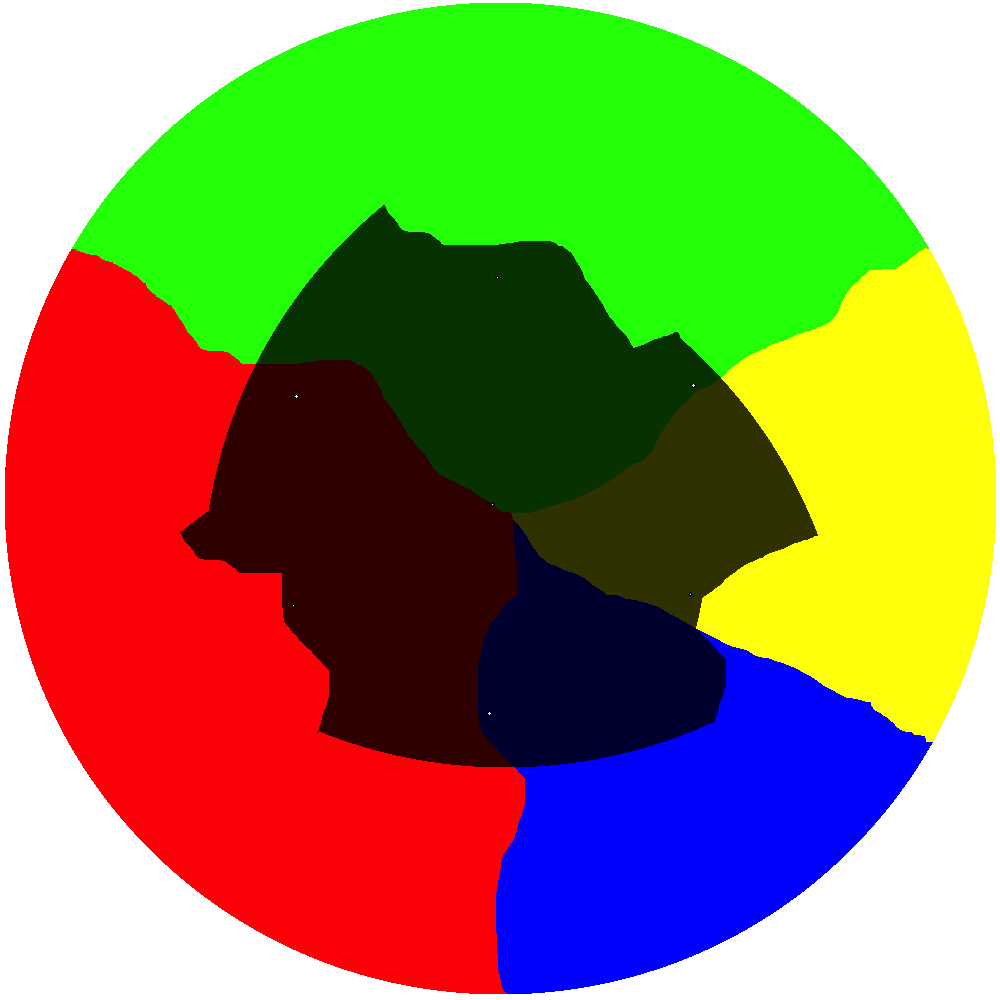} \\
\includegraphics[width=2.5in]{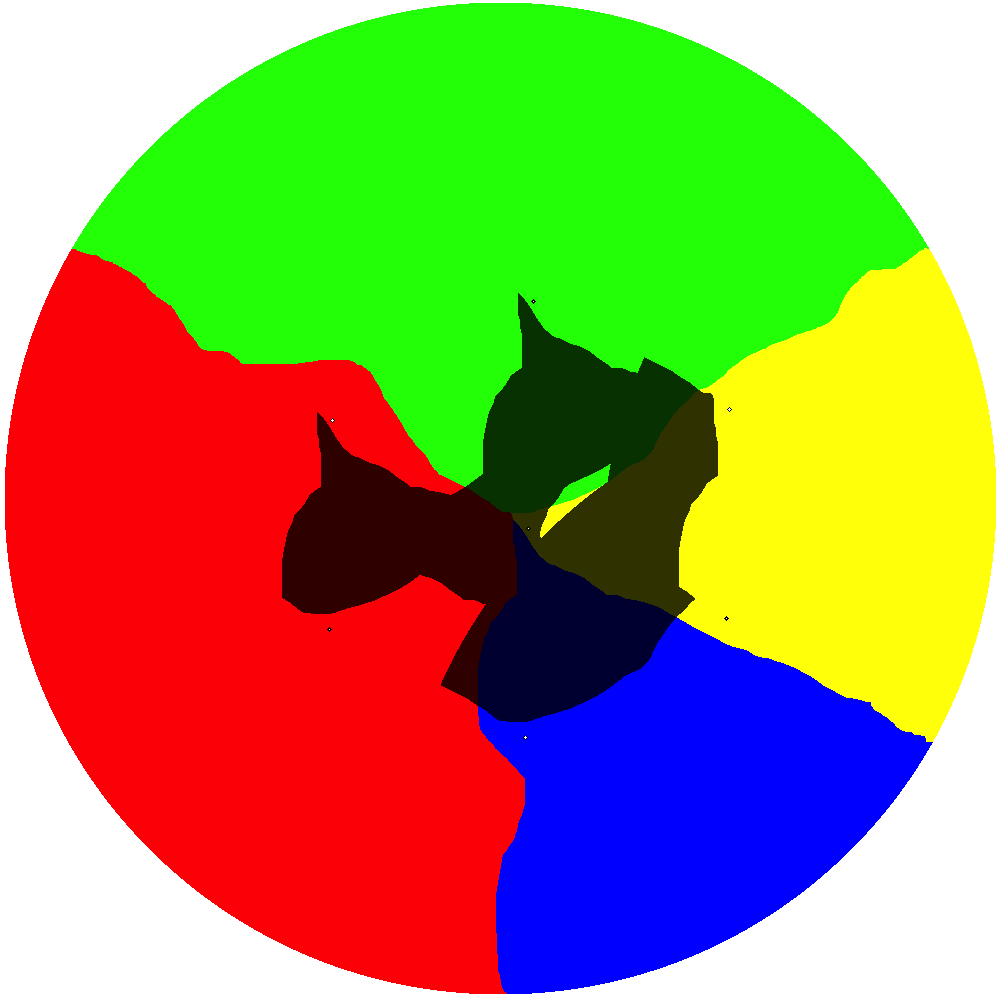} &
\includegraphics[width=2.5in]{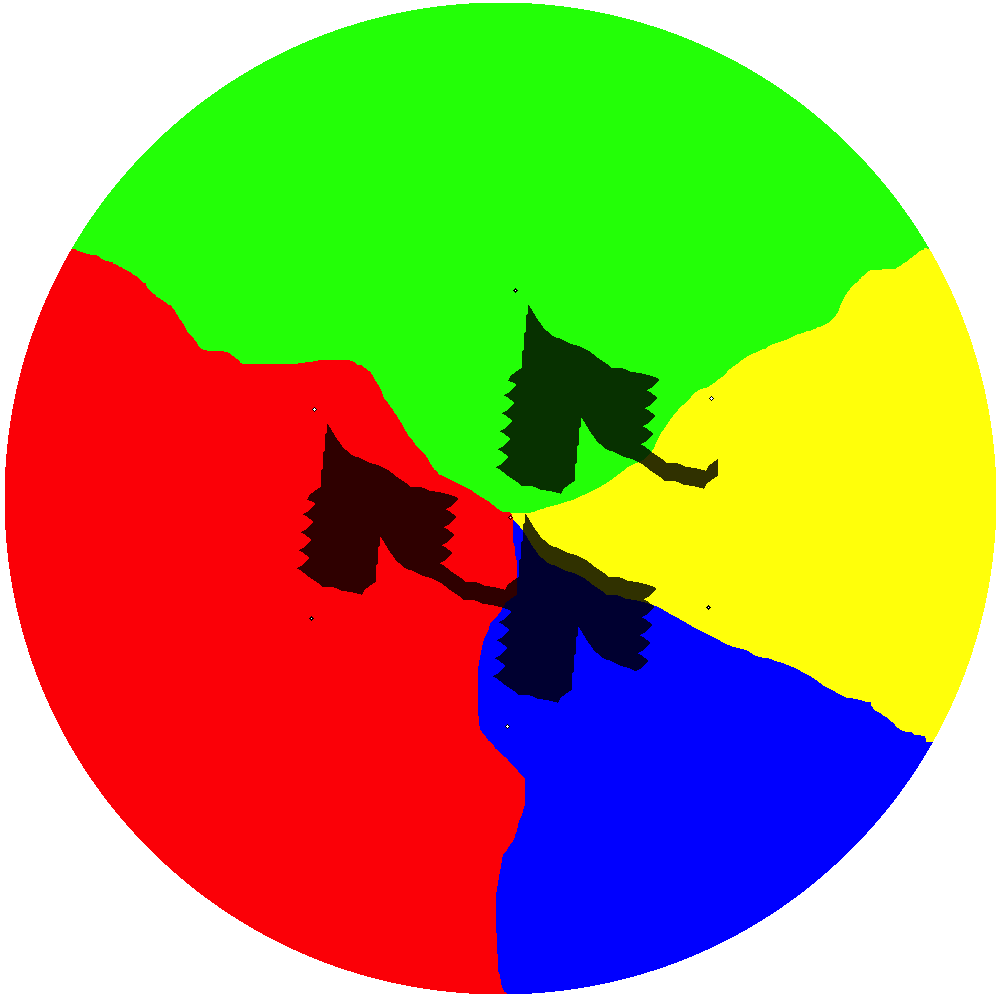}
\end{array}$
\end{center}
\caption{Seemingly finite type attractor for $r = 4$ in $\bbR^2$.\\
The disk~$\Om$ converges to attractor in 28 iterates.
$F^n(\Omega)$ for $n = 1, 2, 5, 28$.}
\label{f:disk-4-fin}
\end{figure}

\subsubsection{Continuity of the attractor\debug{\textbf{Can we say anything smart yet???}}}

Apparently, in general it is not true that $A(T)$ is continuous in the Hausdorff topology. For $m=3, d = 2$, for a certain partitions we numerically observed that in a 1-parameter family of translation vectors, sometimes a large piece of attractor suddenly ``reglues'' from one place to another. On the other hand, these parameter values appear to be discrete, so

\begin{Pbm}
Is it true that at a generic parameter value in a generic family of piecewise translations, $A$ depends continuously on $\mathfrak{v}$? On partition~$\mathfrak{P}$?

\end{Pbm}

\section{Piecewise translations on the $2$-torus: numerics}

As we could not quite deal with piecewise translation maps of $m > d+1$ branches, in this and the following Sections we discuss a number of reductions and simplifications.

For $d = 1$ and $m = 3$, the crucial step to prove that almost every PWT of this kind is of finite type was the following reduction~\cite{Volk2014a}: every piecewise translation map has an interval where the induced map belongs to a special family of piecewise translation maps, \emph{double rotations}, see~\cite{Suzuki2005}, \cite{Bruin2012}. This 3-parameter family of interval maps is given by
\begin{equation}	\label{e:double-rotation}
T_{\alpha, \beta, \delta} (x) := \left\{
           \begin{array}{ll}
             x + \alpha + \beta \mod 1, & \hbox{if $x \le \delta$;} \\
             x + \alpha \mod 1, & \hbox{if $x > \delta$.}
           \end{array}
         \right.
\end{equation}
They can be viewed as piecewise translations of the unit circle~$\bbS = \bbR / \bbZ$ with 2 branches, thus the name.

Though we could not prove a similar reduction theorem for~$d>1$, the idea of the torus rotation factor (recall Subsection~\ref{ss:defs}) inspired us to consider generalizations of double rotations in~$\bbS^d$, $d > 1$.

\subsection{Skew product piecewise translations on the $2$-torus}

Let~$\bbS^2$ be the torus~$\bbR^2 / \bbZ^2$. Let the partition~$\bbS^2 = \Delta_{11} \cup \Delta_{12} \cup \Delta_{21} \cup \Delta_{22}$ be the product of the partitions~$\bbR/\bbZ = \bbS = \delta^x_1 \cup \delta^x_2$ and~$\bbS = \delta^y_1 \cup \delta^y_2$. Consider the class of piecewise translation maps on $\bbS^2$ which are skew products over the circle rotation at the angle~$|\delta^x_2|$:
$$
F(x,y) = (R_{|\delta^x_2|} x, T(x,y)).
$$

Let~$\alpha = |\delta^y_2|$, and $\beta > 0$ some small number.
Define $T(x,y) = T_1(y): = y + \alpha \mod 1$ for $x \in \delta^x_1$. For $x \in \delta^x_2$, let~$T(x,y) = T_2(y) := T_{\alpha, \beta, |\delta_1^y|} (y)$, see~\eqref{e:double-rotation}. So in the fibers over~$\delta^x_1$ we have a rigid rotation, and over~$\delta^x_2$ we have a double rotation~$T_{\alpha, \beta, |\delta_1^y|}$.

This is one of the simplest 2-dimensional piecewise translations of 4 branches, and we believe it is a good model to study. We consider only the nontrivial case when the base rotation is irrational, so that the orbit of every point $x$ in the base is uniformly distributed in~$[0,1]$. The main conjecture suggested by our numerics is that this model is similar to 1-dimensional piecewise translations of 3 branches (see~\cite{Volk2014a}):

\begin{Con}	\label{cj:sk-pr-finite-type}
Generic skew products of this kind are finite type.
\end{Con}


\subsection{Double rotations of a torus}

A somewhat more general setting for piecewise translations of 4 branches on the 2-dimensional torus~$\bbS^2$ is the following. Let $R \subset \bbS^2$ be a rectangle. Let~$\gamma_1, \gamma_2$ be two vectors in $\bbR^2$. Define the piecewise translation map as
$$
F(x) = \left\{
           \begin{array}{ll}
             x + \gamma_1, & \hbox{if $x \in \bbS^2 \setminus R$;} \\
             x + \gamma_2, & \hbox{if $x \in R$.}
           \end{array}
         \right.
$$

\begin{figure}[h]
\begin{center}$
\begin{array}{cc}
\includegraphics[width=2.5in]{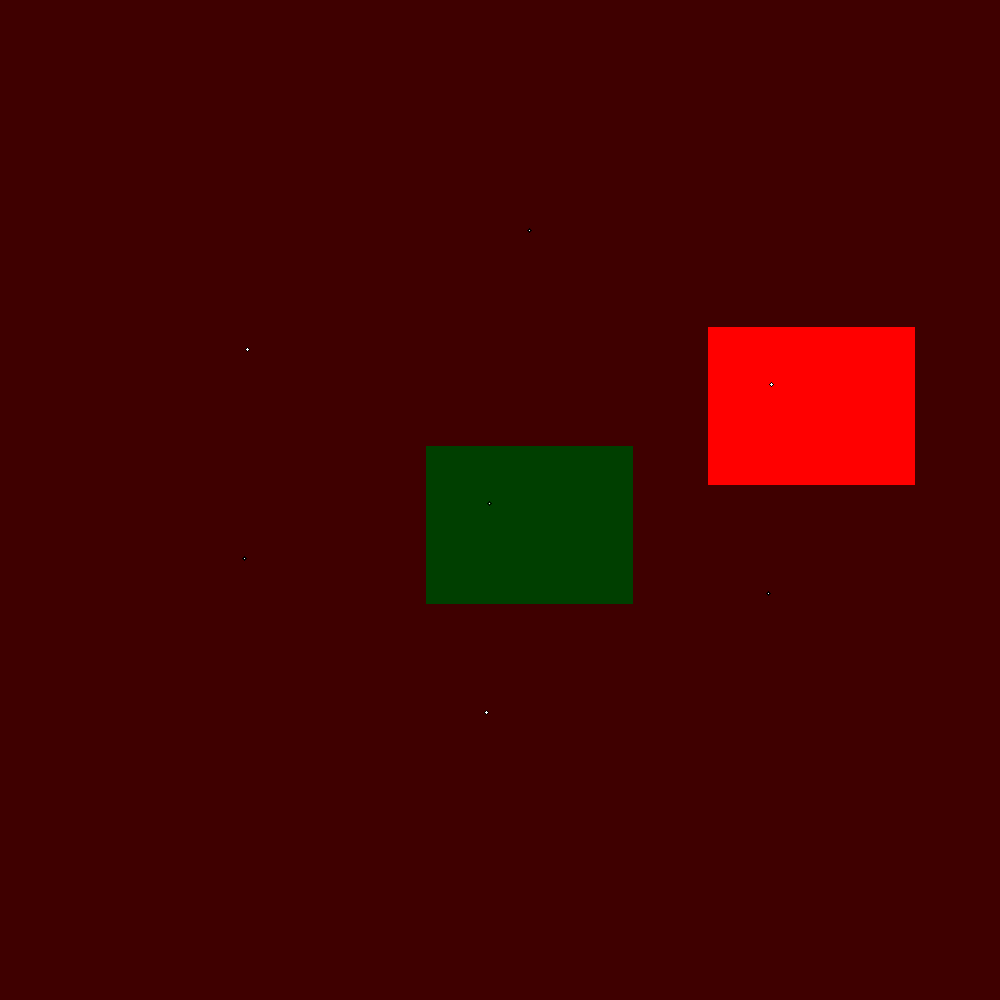} &
\includegraphics[width=2.5in]{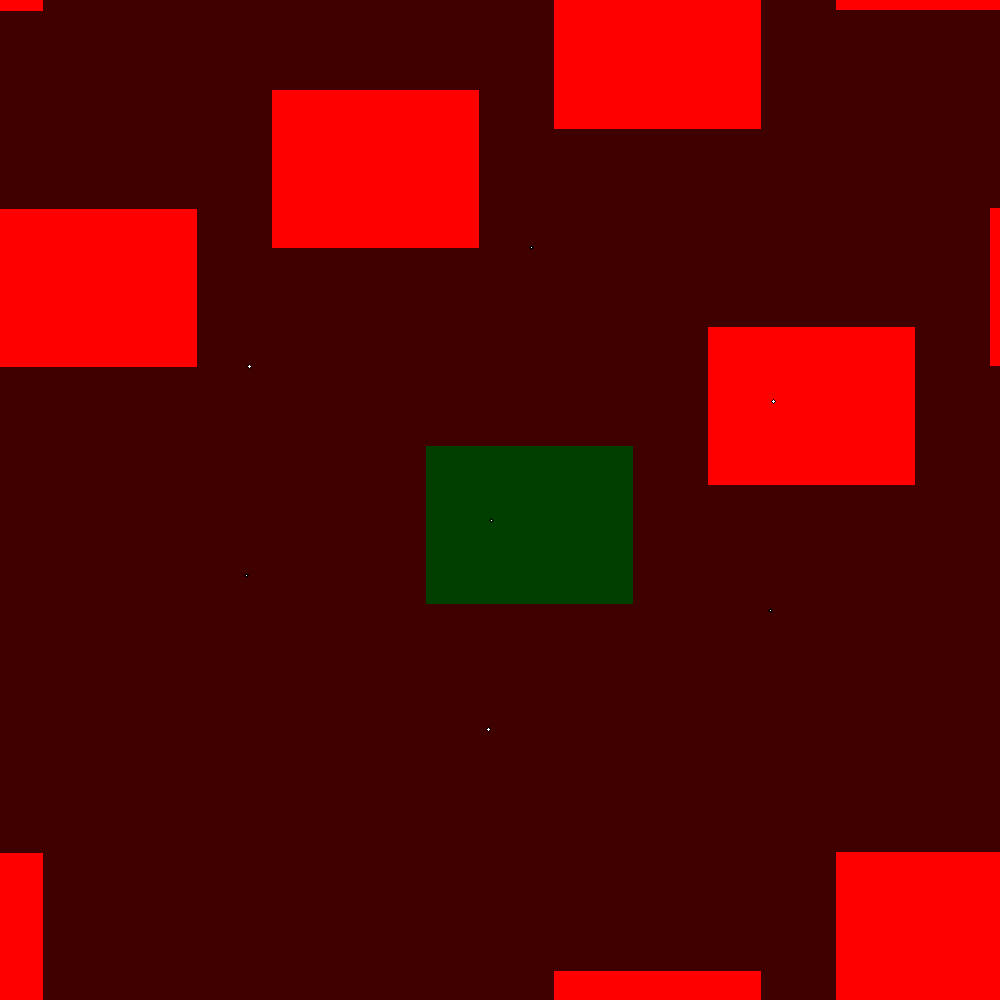} \\
\includegraphics[width=2.5in]{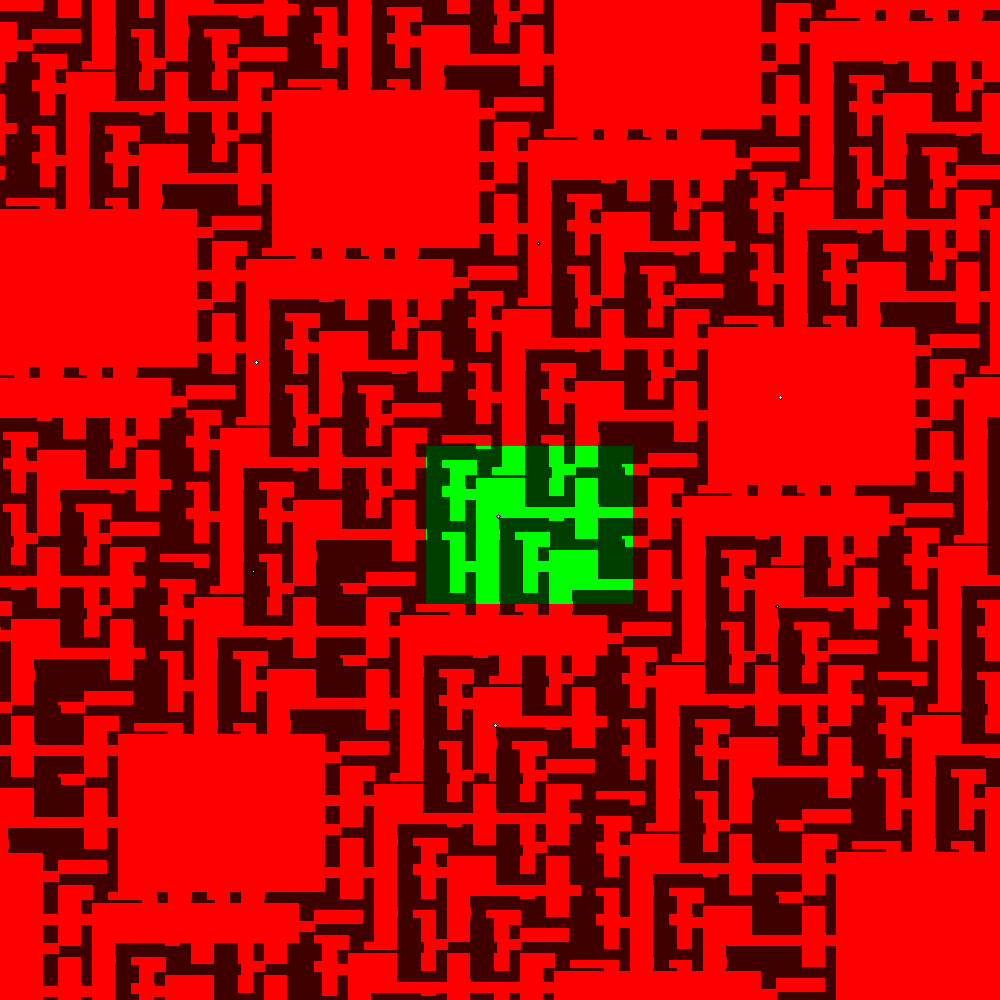} &
\includegraphics[width=2.5in]{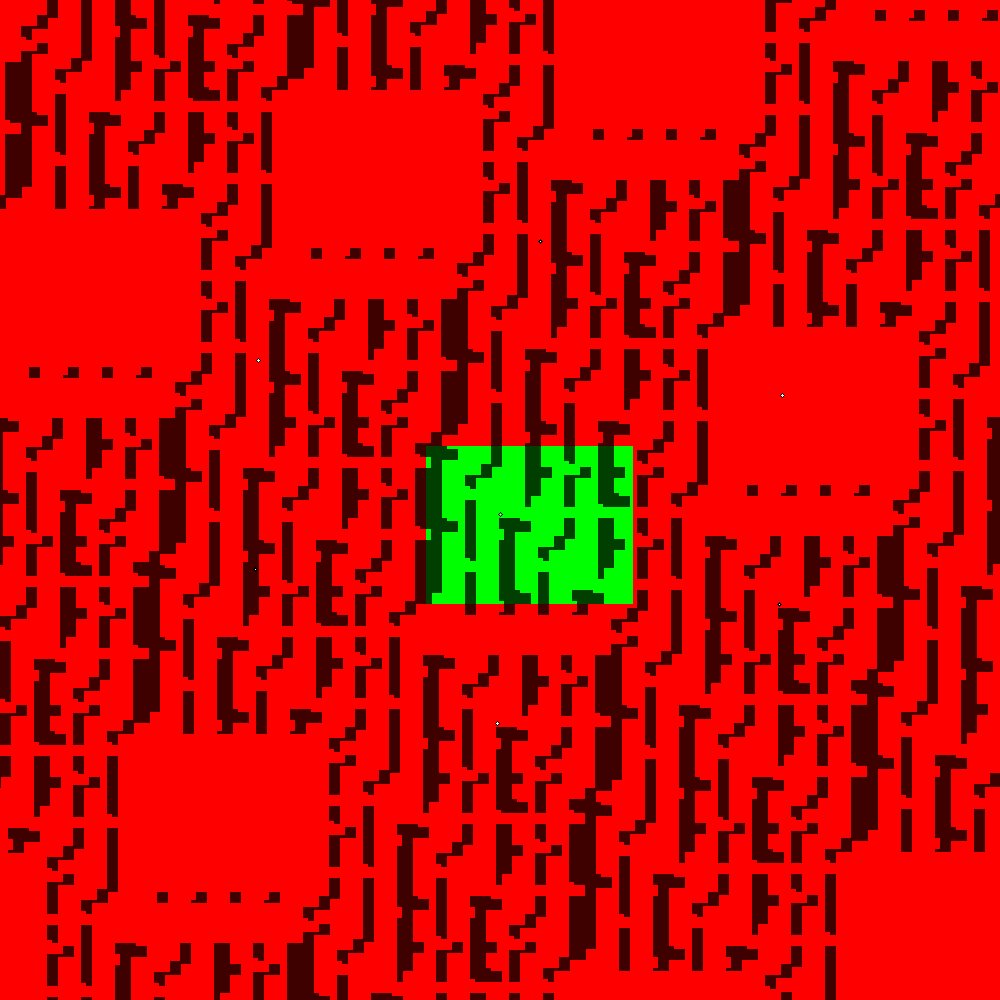}
\end{array}$
\end{center}
\caption{Square on the torus: finite type.\\
Converges to attractor (numerically) in 425 iterates.\\
$F^n (\Omega)$ for $n = 1, 5, 100, 425$.}
\label{f:torus-fin}
\end{figure}

We numerically observed that the behavior of such maps is similar to the one of double rotations~\cite{Bruin2012}.

\begin{Pbm}
\begin{itemize}
\item Prove that for an open, dense, and full measure set of parameters we have finite type attractor.
\item Give an example when attractor is a Cantor set (or has fractal boundary).
\end{itemize}
\end{Pbm}

\begin{figure}[h]
\begin{center}$
\begin{array}{cc}
\includegraphics[width=2.5in]{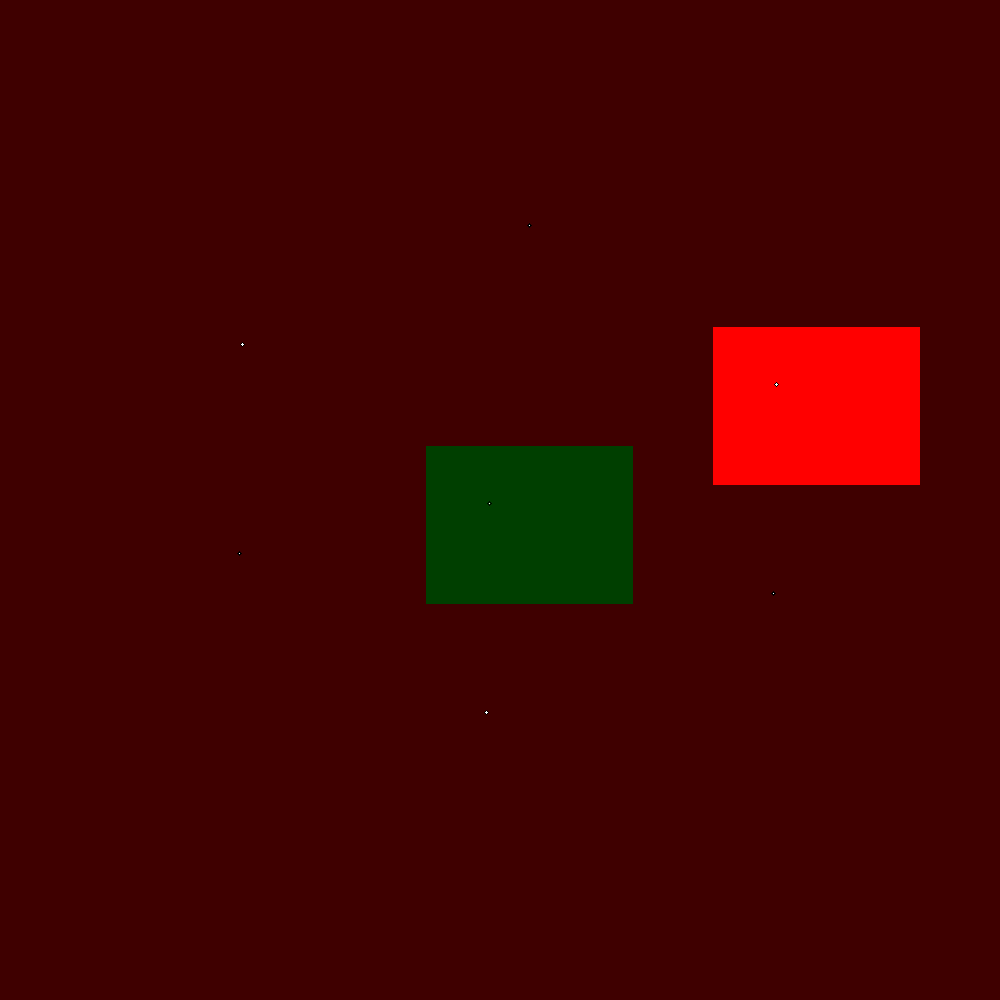} &
\includegraphics[width=2.5in]{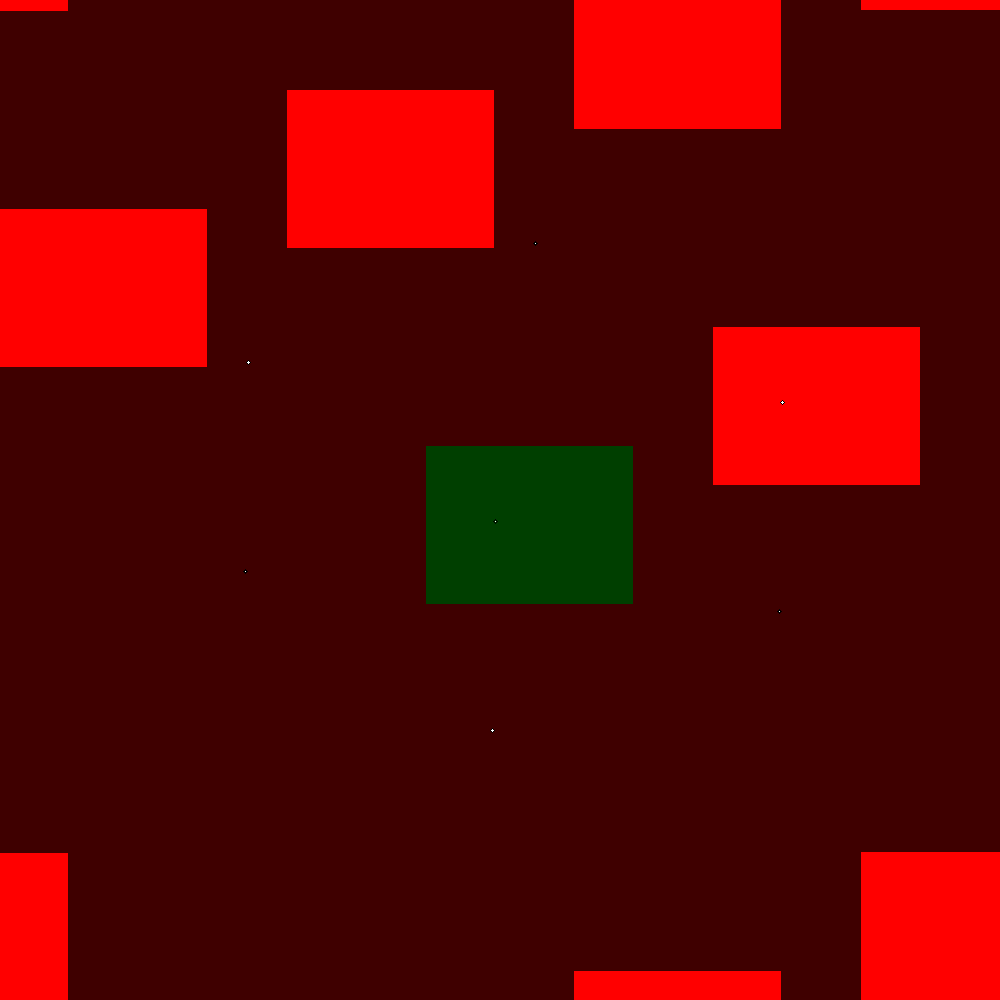} \\
\includegraphics[width=2.5in]{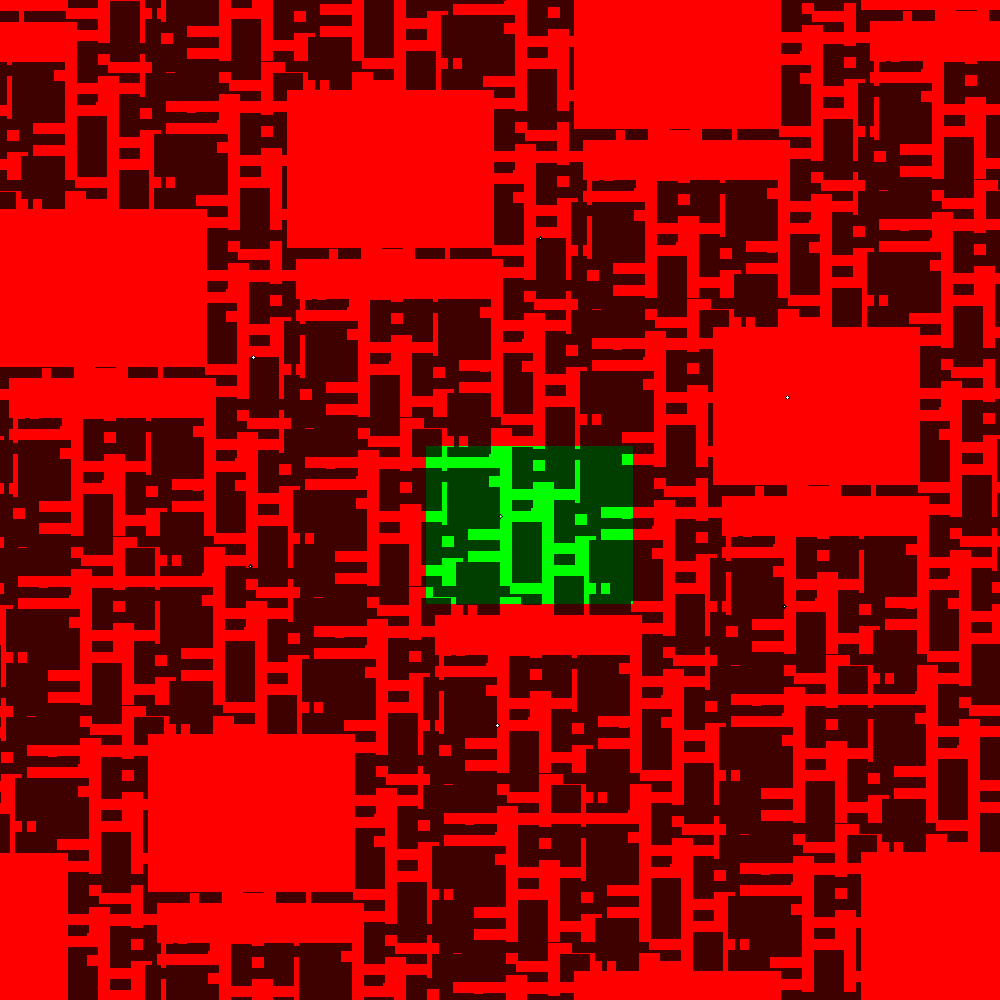} &
\includegraphics[width=2.5in]{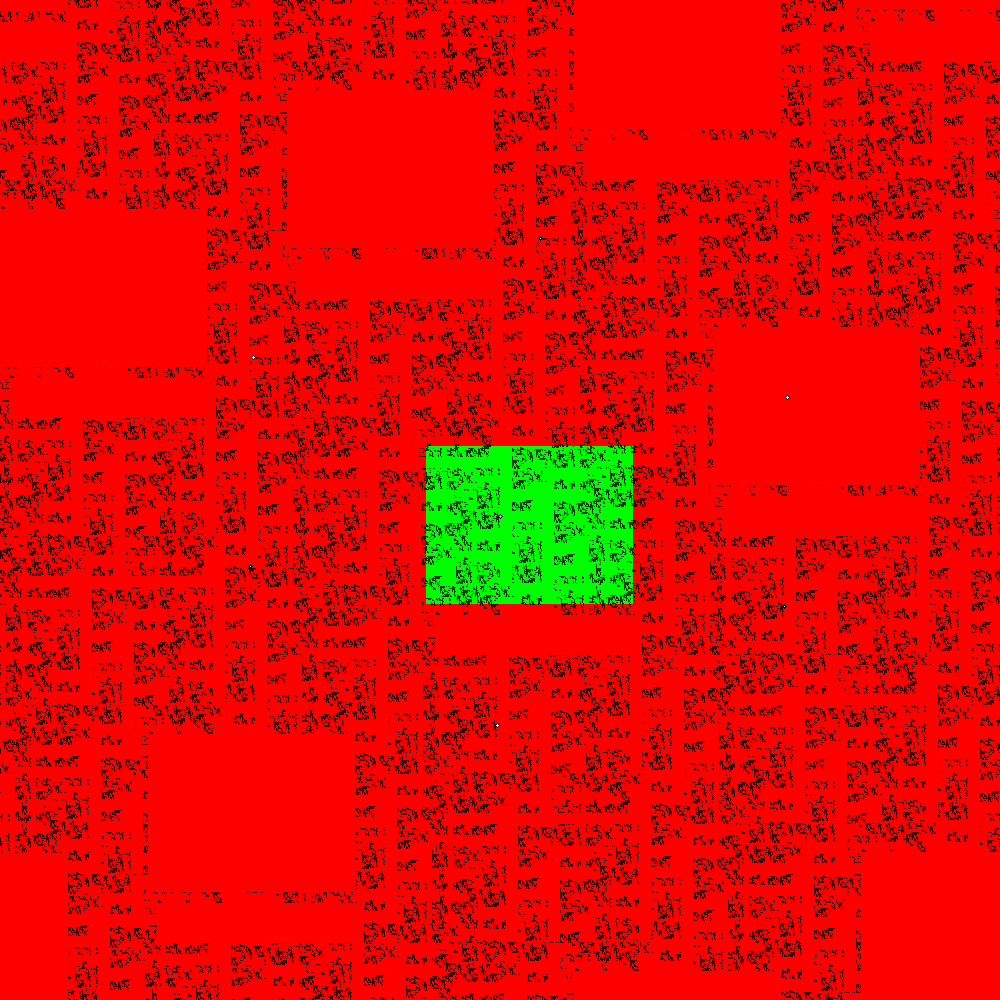}
\end{array}$
\end{center}
\caption{Square on the torus: (seemingly) infinite type.\\
No convergence even in 5000 iterates.
$F^n (\Om)$ for $n = 1, 5, 100, 5000$.}
\label{f:torus-inf}
\end{figure}


\section{Random Double Rotations}

In the previous Section, the partition~$\bbS = \delta^x_1 \cup \delta^x_2$ and the dynamics~$R_{|\delta^x_2|} \colon \bbS \to \bbS$ gave us the fate map~$\mathcal{F} \colon \bbS \ni x \mapsto \om \in \{1,2\}^\bbN$. As $|\delta^x_2| \notin \bbQ$, for any $x \in \bbS$ the image $\mathcal{F}(x)$ is a aperiodic sequence, but it is still deterministic. It is reasonable to consider what happens in the totally random base dynamics.


So,
we take the sequence of Bernoulli~$(p, 1-p)$ independent random variables~$i_n$, and let
$$
F_n = T_{i_1} \circ \dots \circ T_{i_n},
$$
with the same~$T_1, T_2$ as before.

\begin{figure}[h]
\centering
\includegraphics[width=420pt]{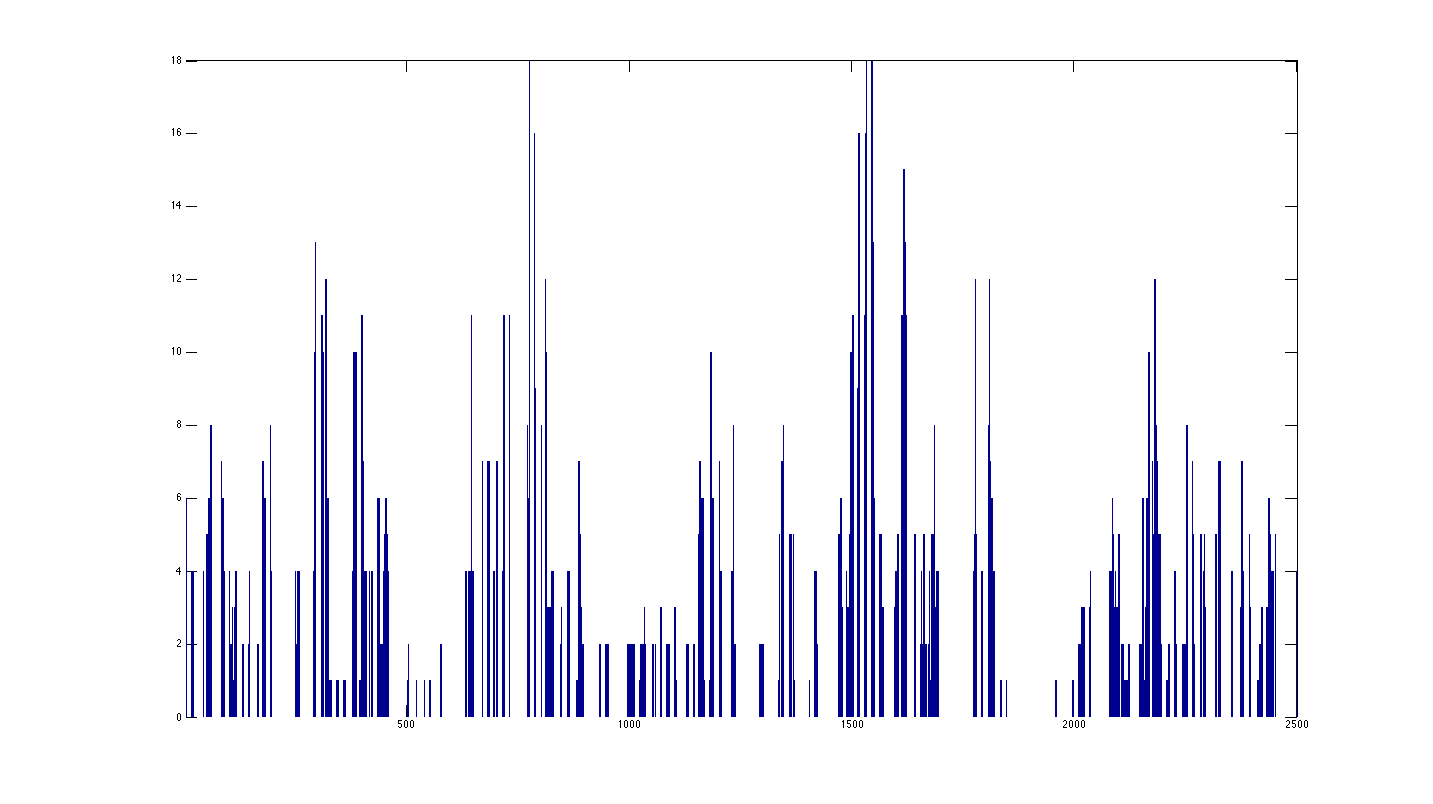}\\
\caption{Random double rotations.\\ Attractor measure after 2000 iterates.}
\label{f:random-double-rotations}
\end{figure}

\begin{Thm} \label{t:leb-to-0}
  Almost surely,
$$
\lim_{n\to\infty} \Leb(F_n \bbS) = 0.
$$
\end{Thm}

This theorem follows from

\begin{Lem}\label{l:arc-itinerary}
For any arc~$A \subset \bbS$, there exists a finite itinerary~$J = j_1, \dots, j_k$ such that~$T_J \bbS := (T_{j_1} \circ \dots \circ T_{j_k}) \bbS \subset A$.
\end{Lem}

\begin{proof}
First, let us show that there exists a finite itinerary~$J = j_1, \dots, j_k$ such that~$T_J = T_{j_1} \circ \dots \circ T_{j_k}$ sends the circle~$\bbS$ into some arc of length less than~$\beta$.
About~$T_2$, note that
\begin{itemize}
\item the set $T_2\bbS$ has a gap~$[\alpha, \alpha+\beta]$ of length~$\beta$;
\item the arc $[\delta_1^y; \delta_1^y - \beta]$ of length~$1-\beta$ is mapped 1-1 onto its image.
\end{itemize}
Now take $J$ of the form
\begin{equation*}
J = 2\underbrace{1\dots 1}_{m_l \text{ times}} 2 \underbrace{1\dots 1}_{m_{l-1} \text{ times}}   2 \dots 2 \underbrace{1\dots 1}_{m_1 \text{ times}} 2.
\end{equation*}
After the first $T_2$, the image of the circle $\bbS$ has an arc gap of size $\beta$. Put in sufficiently many ($m_1$) rotations~$T_1$ so that after $T_1^{m_1} \circ T_2$ the gap
\begin{enumerate}
\item\label{i:1-1} has empty intersection with the non-1-1 arc~$[\delta_1^y - \beta; \delta_1^y]$ of length~$\beta$;
\item\label{i:1-2} overlaps with $0$.
\end{enumerate}
After the next~$T_2$ the image of~$\bbS$ has an arc gap of length~$2\beta$. The procedure can be continued provided we can satisfy the above conditions~\ref{i:1-1}, \ref{i:1-2}. When this eventually becomes impossible, the image $T_J \bbS$ has an arc gap of length~$1-\beta < L \le 1$, which means $T_J$ sends $\bbS$ into an arc of length less than~$\beta$.

Now let us show that for an arbitrary arc~$K$, $|K| < \beta$, there exists a finite itinerary which sends $K$ into some arc of length less than~$\frac34 |K|$. Attach two markers~$\mm_1, \mm_2$ at the $\frac14$ and $\frac34$ of $K$ (so that $K$ is within~$\frac14|K|$-neighborhood of the set~$\{\mm_1,\mm_2\}$), and let us track their orbits.

For this, we construct an itinerary of the form
\begin{equation*}
J = 2\underbrace{1\dots 1}_{m_l \text{ times}} 2 \underbrace{1\dots 1}_{m_{l-1} \text{ times}}   2 \dots 2 \underbrace{1\dots 1}_{m_1 \text{ times}}
\end{equation*}
such that along this itinerary,
\begin{enumerate}
\item the image of $K$ is always either a single arc or two arcs, in the latter case each contains exactly one of the markers;
\item\label{i:2-2} the image of $K$ is always within~$\frac14|K|$-neighborhood of the set~$\{F_n(\mm_1),F_n(\mm_2)\}$;
\item the difference~$F_n(\mm_2) - F_n(\mm_1)$ either does not change or changes by $+\beta$, moreover, the latter happens infinitely many times.
\end{enumerate}

\emph{Phase 1.} Using the irrational rotations~$T_1$, we can position the image of $K$ so that
\begin{itemize}
\item if the image of $K$ is a single arc, which implies $F_n(\mm_2) - F_n(\mm_1) < \frac12|K|$, then $0$ is $\frac12(\frac12|K| - (F_n(\mm_2) - F_n(\mm_1)))$-close to the midpoint of $[F_n(\mm_2),F_n(\mm_1)]$;
\item if the image of $K$ is two arcs, then $0$ is in their complement.
\end{itemize}
Note that the difference~$F_n(\mm_2) - F_n(\mm_1)$ does not change with $T_1$'s.

\emph{Phase 2.} Use the double rotation~$T_2$. Note that the difference $F_n(\mm_2) - F_n(\mm_1)$ just increased by $\beta$. If the image of $K$ was two arcs before this, the condition~\ref{i:2-2} holds trivially by induction. If the image of $K$ was a single arc, then after the double rotation the condition~\ref{i:2-2} holds by the previous positioning of $F_n(K)$ with respect to $0$.

Repeat \emph{Phase 1.}

Thus we construct an itinerary which makes $\mm_2$ to do infinite irrational rotation at the angle~$\beta$ with respect to $\mm_1$, so they eventually get $\frac14|K|$-close. Then because of condition~\ref{i:2-2}, the whole image $F_n(K)$ is contained within an arc of length~$\frac34|K|$.

Finally, repeating the process, for any $m \in \bbN$ we can send $\bbS$ into an arc of length less than~$\left(\frac34\right)^m\cdot|K|$. Taking $m$ sufficiently large and subsequently using many rotations, we can also send the result into~$A$.
\end{proof}

\begin{proof}[Proof of Theorem~\ref{t:leb-to-0}]
By Lemma~\ref{l:arc-itinerary}, for an arbitrary small interval~$\dd$ there exists a finite itinerary which sends~$\bbS$ into~$\dd$. This itinerary is met in almost every infinite itinerary with the probability one (actually, even with a positive frequency). As the Lebesgue measure of $F_n\bbS$ is monotone nonincreasing as $n\to\infty$, this means that the limit is almost surely less than any positive number, thus being equal to zero.
\end{proof}

So for the random iterations of a rigid rotation and a double rotation, the attractor is almost surely a closed set of zero Lebesgue measure. This can be viewed as a nonhyperbolic syncronization phenomenon: the initial conditions bunch up and are eventually concentrated on some tiny subset of the phase space. A similar behavior was observed for random iterations of generic diffeomorphisms of a circle in~\cite{Antonov1984},~\cite{Kleptsyn2004}, and of an interval in~\cite{Kleptsyn2014a}.

\begin{Pbm}
In the above setting, prove that the attractor is almost surely a Cantor set (i.e. there are no isolated points). What is its Hausdorff dimension?
\end{Pbm}


\begin{Rem}
\begin{itemize}
\item By the Kolmogorov's zero-one law, either the attractor is almost surely a Cantor set, or it almost surely is not.
\item Same about attractor being a finite set of points.
\item Being just a single point is not a Kolmogorov's zero-one property.
\item For a similar reason (ergodicity of Bernoulli shift), the dimension of the attractor is almost surely some constant depending on the maps~$T_1, T_2$.
\end{itemize}
\end{Rem}

\begin{Pbm}
Averaging the Lebesgue measure over iterates and all possible fates gives some measure supported on the attractor. What could be said about this measure?
\end{Pbm}


\section*{Acknowledgements}

The authors thanks A.~Gorodetski and UC Irvine for support and hospitality during the initial phase of the project, T.~Nowicki, J.~Schmeling, M.~Welling, W.~Yessen for insightful discussions and S.~Truong for engineering references.

\bibliographystyle{plain}
\bibliography{dynsys-utf8}

\def\cprime{$'$}
\begin{thebibliography}{10}

\bibitem{Adler2012}
R.~Adler, T.~Nowicki, G.~{\'S}wirszcz, C.~Tresser, and S.~Winograd.
\newblock Error diffusion on simplices: Invariant regions, tessellations and
  acuteness.
\newblock {\em Preprint}, 2012.

\bibitem{Adler2005}
R.~L. Adler, B.~Kitchens, M.~Martens, C.~Pugh, M.~Shub, and C.~Tresser.
\newblock Convex dynamics and applications.
\newblock {\em Ergodic Theory Dynam. Systems}, 25(2):321--352, 2005.

\bibitem{Adler2010}
R.~L. Adler, T.~Nowicki, G.~{\'S}wirszcz, and C.~Tresser.
\newblock Convex dynamics with constant input.
\newblock {\em Ergodic Theory Dynam. Systems}, 30(4):957--972, 2010.

\bibitem{akiyama2015pisot}
Shigeki Akiyama, Marcy Barge, Val{\'e}rie Berth{\'e}, J-Y Lee, and Anne Siegel.
\newblock On the pisot substitution conjecture.
\newblock In {\em Mathematics of Aperiodic Order}, pages 33--72. Springer,
  2015.

\bibitem{Antonov1984}
V.~A. Antonov.
\newblock Modeling of processes of cyclic evolution type. {S}ynchronization by
  a random signal.
\newblock {\em Vestnik Leningrad. Univ. Mat. Mekh. Astronom.}, 2:67---76, 1984.

\bibitem{ashwin2001dynamics}
Peter Ashwin, Jonathan~HB Deane, and X-C Fu.
\newblock Dynamics of a bandpass sigma-delta modulator as a piecewise isometry.
\newblock In {\em Circuits and Systems, 2001. ISCAS 2001. The 2001 IEEE
  International Symposium on}, volume~3, pages 811--814. IEEE, 2001.

\bibitem{avila2007weak}
Artur Avila and Giovanni Forni.
\newblock Weak mixing for interval exchange transformations and translation
  flows.
\newblock {\em Annals of Mathematics}, pages 637--664, 2007.

\bibitem{Boshernitzan1995}
Michael Boshernitzan and Isaac Kornfeld.
\newblock Interval translation mappings.
\newblock {\em Ergodic Theory Dynam. Systems}, 15(5):821--832, 1995.

\bibitem{Bruin2012}
Henk Bruin and Gregory Clack.
\newblock Inducing and unique ergodicity of double rotations.
\newblock {\em Discrete and Continuous Dynamical Systems - Series A},
  32(12):4133--4147, 2012.

\bibitem{chua1988chaos}
Leon~O Chua and Tao Lin.
\newblock Chaos in digital filters.
\newblock {\em IEEE Transactions on Circuits and Systems}, 35(6):648--658,
  1988.

\bibitem{deane2004buck}
JHB Deane.
\newblock The buck converter planar piecewise isometry.
\newblock {\em Preprint}, pages 1--4, 2004.

\bibitem{deane2002global}
Jonathan~HB Deane.
\newblock Global attraction in the sigma--delta modulator piecewise isometry.
\newblock {\em Dynamical Systems}, 17(4):377--388, 2002.

\bibitem{deane2006piecewise}
Jonathan~HB Deane.
\newblock Piecewise isometries: applications in engineering.
\newblock {\em Meccanica}, 41(3):241--252, 2006.

\bibitem{feely2000nonlinear}
Orla Feely, Daniele Fournier-Prunaret, Ina Taralova-Roux, and David Fitzgerald.
\newblock Nonlinear dynamics of bandpass sigma-delta modulation—an
  investigation by means of the critical lines tool.
\newblock {\em International Journal of Bifurcation and Chaos},
  10(02):307--323, 2000.

\bibitem{Kleptsyn2014a}
V.~Kleptsyn and D.~Volk.
\newblock Physical measures for nonlinear random walks on interval.
\newblock {\em Moscow Mathematical Journal}, 14(2):339--365, 2014.

\bibitem{Kleptsyn2004}
V.~A. Kleptsyn and M.~B. Nalskii.
\newblock Contraction of orbits in random dynamical systems on the circle.
\newblock {\em Functional Analysis and Its Applications}, 38(4):267--282,
  October 2004.

\bibitem{masur1982interval}
Howard Masur.
\newblock Interval exchange transformations and measured foliations.
\newblock {\em Annals of Mathematics}, 115(1):169--200, 1982.

\bibitem{Schmeling2000}
{J\"org} Schmeling and Serge Troubetzkoy.
\newblock Interval translation mappings.
\newblock In {\em Dynamical systems ({L}uminy-{M}arseille, 1998)}, pages
  291--302. World Sci. Publ., River Edge, NJ, 2000.

\bibitem{Suzuki2004}
H.~Suzuki, K.~Aihara, and T.~Okamoto.
\newblock Complex behavior of a simple partial-discharge model.
\newblock {\em Europhysics Letters}, 66(1):28--24, 2004.

\bibitem{Suzuki2005}
Hideyuki Suzuki, Shunji Ito, and Kazuyuki Aihara.
\newblock Double rotations.
\newblock {\em Discrete Contin. Dyn. Syst.}, 13(2):515--532, 2005.

\bibitem{veech1982gauss}
William~A Veech.
\newblock Gauss measures for transformations on the space of interval exchange
  maps.
\newblock {\em Annals of Mathematics}, 115(2):201--242, 1982.

\bibitem{viana2006ergodic}
Marcelo Viana.
\newblock Ergodic theory of interval exchange maps.
\newblock {\em Revista Matem{\'a}tica Complutense}, 19(1):7--100, 2006.

\bibitem{Volk2014a}
Denis Volk.
\newblock Almost every interval translation map of three intervals is finite
  type.
\newblock {\em Discrete Contin. Dyn. Syst.}, 34(5):2307--2314, 2014.

\bibitem{Welling2009}
Max Welling.
\newblock Herding dynamic weights for partially observed random field models.
\newblock In {\em In Proc. of the Conf. on Uncertainty in Artificial
  Intelligence}, 2009.

\end{thebibliography}

\end{document}